\documentclass[11pt]{article}
\usepackage{amssymb,amsmath,bbm,latexsym,amscd}
\usepackage{stmaryrd} % pour les \llbracket et \rrbracket
\usepackage[mathscr]{euscript} % pour des belles lettres
\usepackage{relsize} % pour utiliser \mathlarger{\mathlarger{}} sur des lettres de \mathscr

\usepackage{color}
\definecolor{orange}{rgb}{1, 0.3, 0.1}

\newtheorem{theorem}[equation]{Theorem}
\newtheorem{corollary}[equation]{Corollary}
\newtheorem{lemma}[equation]{Lemma}
\newtheorem{proposition}[equation]{Proposition}
\newtheorem{definition}[equation]{Definition}
\newtheorem{remark}[equation]{Remark}
\numberwithin{equation}{section}

\newcommand{\bb}{\mathbbm{b}}
\newcommand{\ot}{\otimes}

\newcommand{\Cx}{\mathbb{C}}

\newcommand{\fg}{{\mathfrak g}}
\newcommand{\fh}{{\mathfrak h}}

\newcommand{\Rep}{\hbox{Rep}\,}
\DeclareMathOperator{\Max}{Max} % extensions

\newcommand{\ga}{\alpha}

\newcommand{\proof}{{\bf Proof\ \ }}
\newcommand{\qed}{\hfill $\Box$}

\newcommand{\LL}{{\mathcal L}}

\newcommand{\F}{{\mathcal F}}

\newcommand{\B}{{\mathcal B}}
\newcommand{\End}{\hbox{End}}

\newcommand{\cL}{{\mathcal L}}
\newcommand{\R}{{\mathcal R}}

\newcommand{\Aut}{{\rm Aut}}

\newcommand{\Hom}{{\rm Hom}}

\newcommand{\Span}{{\rm Span}}

\newcommand\Spec{\text{\rm Spec}\,}

\newcommand{\ev}{\mathrm{ev}}
\DeclareMathOperator{\supp}{supp} % support
\DeclareMathOperator{\Ext}{Ext} % extensions
\DeclareMathOperator{\Der}{Der} % derivations
\DeclareMathOperator{\IDer}{IDer} % derivations interieures
\DeclareMathOperator{\HHH}{H} % H de cohomologie
\newcommand{\cohomo}[2]{\HHH^1(#1 \; ; #2)} % premiere cohomologie
\newcommand{\gcohomo}[2]{\HHH^1\big(#1 \; ; #2\big)} % premiere cohomologie (grosses parentheses)
\DeclareMathOperator{\Bloc}{Ext-block} % bloc d'extension
\newcommand{\secr}{\mathlarger{\mathscr{F}}} % section fibré irreps
\newcommand{\carspec}{\mathlarger{\mathscr{SC}}} % caractère spectraux
\DeclareMathOperator{\nev}{comp-ev} % non-d'évaluation... (un module)

\title{Extensions of modules for twisted current algebras}

\date{}

\author{Jean Auger$\hbox{}^{1* }$ and Michael Lau$\hbox{}^{2 }$\thanks{The authors gratefully acknowledge funding received from the Fonds de Recherche du Qu\'ebec-Nature et Technologies (J.~Auger) and the Natural Sciences and
Engineering Research Council of Canada (M.~Lau).} \vspace{0.3cm}
\vspace{0.3cm}\\
$\hbox{\ \,}^1${\small University of Alberta}\\ {\small Department of
Mathematical and Statistical Sciences}\\{\small Edmonton, AB,
Canada T6G 2G1}\\
\\
$\hbox{}^{2 }${\small Universit\'e Laval}\\
{\small D\'epartement de math\'ematiques et de statistique}\\ {\small Qu\'ebec, QC, Canada G1V 0A6}\vspace{0.1cm}}

\begin{document}
\maketitle

\begin{small}
\noindent {\bf Abstract.} Twisted current algebras are fixed point subalgebras of current algebras under a finite group action.  Special cases include equivariant map algebras and twisted forms of current algebras.  Their finite-dimensional simple modules fall into two categories, those which factor through an evaluation map and those which do not.  We show that there are no nontrivial extensions between finite-dimensional simple evaluation and non-evaluation modules.  We then compute extensions between any pair of finite-dimensional simple modules for twisted current algebras, and use this information to determine the block decomposition for the category.  In the special case of twisted forms, this decomposition can be described in terms of maps to the fundamental group of the underlying root system.

%We compute the extensions between any pair of simple modules in the category of finite-dimensional modules over a twisted current algebra, and then use this information to determine the block decomposition of the category. We illustrate our results with an application to twisted forms of current algebras.

%We use evaluation representations to give a complete classification of the finite-dimensional simple modules of twisted current algebras.  This generalizes and unifies recent work on multiloop algebras, current algebras, equivariant map algebras, and twisted forms.

\bigskip

\noindent {\bf Keywords:} current algebras, finite-dimensional modules, extensions of modules, block decomposition, twisted current algebras

\bigskip

\noindent
{\bf MSC2010:} primary 17B55; secondary 18G15, 17B10, 17B65

\end{small}
\maketitle

%%%
\section{Introduction}
Let $\fg$ be a finite-dimensional reductive Lie algebra over an algebraically closed field $k$ of characteristic zero, and suppose that $X$ is an affine $k$-scheme of finite type with coordinate algebra $S$.  The vector space $\fg$ is naturally an affine $k$-scheme, and the algebra of morphisms from $X$ to $\fg$ may be identified with the {\em current Lie algebra} $\fg\ot_k S$, under the pointwise Lie bracket $[x\ot r,y\ot s]=[x,y]\ot rs$ for all $x,y\in\fg$ and $r,s\in S$.

Let $\Gamma$ be a finite group acting by $k$-algebra automorphisms on $S$, and let $u$ be a crossed homomorphism from $\Gamma$ to the group $\Aut_{S-Lie}(\fg\ot_k S)$ of $S$-Lie automorphisms of $\fg\ot_k S$.  The crossed homomorphism $\gamma\mapsto u_\gamma$ then induces a group action of $\Gamma$ on the current algebra $\fg\ot_k S$ by $k$-Lie algebra automorphisms:
$${}^{\gamma}(x\ot s):=u_\gamma(x\ot{}^{\gamma}s),$$
and the $k$-Lie subalgebra   $(\fg\ot_k S)^\Gamma=\{z\in\fg\ot S\ :\ {}^\gamma z=z\hbox{\ for all\ } \gamma\in\Gamma\}$ of $\Gamma$-invariants is called the {\em twisted current algebra} associated with the action of $\Gamma$ on $\fg\ot_k S$.  In the case where $\fg$ is simple and $S$ is reduced, any group action $\Gamma\rightarrow\hbox{Aut}_{k-Lie}(\fg\ot_k S)$ is of this form, and the twisted current algebras are precisely the invariant subalgebras of $\fg\ot_k S$ under finite group actions \cite[Theorem 2.2]{twcurr}.

%In the case where $\fg$ is simple and $S$ is reduced, the action of $\Gamma$ on $\fg\ot S$ induces a group action on $S$ by $k$-algebra automorphisms, and a crossed homomorphism $u:\ \Gamma\rightarrow\Aut_{S-Lie}(\fg\ot S)$, such that
%\begin{equation}\label{rem:splitting}
%{}^\gamma(x\ot s)=u_\gamma(x\ot {}^\gamma s),
%\end{equation}
%for all $\gamma\in\Gamma$, $x\in\fg$, and $s\in S$ \cite[Proposition 2.2]{twcurr}.

We are especially interested in the case where $\fg$ is simple and $S$ is a Galois extension of the subalgebra $S^\Gamma$ of $\Gamma$-invariants.  Descent theory studies these twisted current algebras, which are exactly the {\em $S/S^\Gamma$-twisted forms} of the current algebra $\fg\ot_k S^\Gamma$.  Modules for these algebras will be studied in detail in Section 5.  Another important class of examples is the case where the cocycle $u$ defines a group action of the form ${}^\gamma(x\ot s)={}^\gamma x\ot {}^\gamma s$ for all $\gamma\in\Gamma$ and $x\ot s\in\fg\ot S$.  Such twisted current algebras are called {\em equivariant map algebras}.  For further examples and details on twisted current algebras, see \cite{twcurr}.

The finite-dimensional simple modules for twisted current algebras were classified in \cite{twcurr} as {\em evaluation modules} or tensor products of evaluation modules with 1-dimensional {\em non-evaluation} modules. Evaluation modules are those whose action factors through the image of an evaluation homomorphism
$$\ev_{\mathbf{M}}:\ \cL\hookrightarrow\fg\ot_k S\rightarrow \left(\fg\ot_k S/M_1\right)\oplus\cdots\oplus\left(\fg\ot_k S/M_r\right)$$
for some family of maximal ideals $\mathbf{M}=\{M_1,\ldots,M_r\}\subseteq\Max S$.  As the image of each of the ``evaluations'' $\cL\rightarrow \fg\ot S/M_i$ is isomorphic to a reductive subalgebra of the finite-dimensional $k$-Lie algebra $\fg$, this reduces the classification of evaluation representations to the well known description of simple modules for reductive Lie algebras and the somewhat trickier problem of classification up to isomorphism.  Non-evaluation modules first appeared in the context of the Onsager algebra \cite{DR}; they are representations which do not factor through any evaluation homomorphism.

%whose kernel does not contain any finite intersection of kernels of evaluation maps $\ev_M$, for maximal ideals $M\in\Max S$.

The next step in understanding the category $\LL$-$\mathbf{mod}$ of finite-dimensional modules for twisted current algebras, and the main goal of the present article, is to calculate the extensions between simple modules and determine the blocks of $\LL$-$\mathbf{mod}$.  Such a classification was recently obtained in the case of equivariant map algebras \cite{NS}, modulo an assumption that certain ext-groups vanish.

\bigskip

After a brief review of the classification of finite-dimensional simple modules in Section 2, we address the open question of calculating extensions between evaluation and non-evaluation modules in Section 3.  The block decomposition for equivariant map algebras was previously obtained only under the key assumption that $\Ext_\cL^1(V,k_\lambda)=0$ for simple evaluation modules $V$ and $1$-dimensional non-evaluation modules $k_\lambda$.  Such algebras $\cL$ are said to be {\em extension-local}.  %In Section 2, we review various elementary facts pertinent to the calculation of the space $\Ext_\cL^1(V,W)$ of extensions between $\cL$-modules $V$ and $W$.  In Section 3, we rapidly recall the notation and main theorems of \cite{twcurr} needed for the classification of ext-blocks.  Section 4 considers extensions between evaluation modules.  We restate several of the main results of \cite{NS}, extending them to the strictly more general context of twisted current algebras.  Most of the proofs in this section are straightforward generalizations of the equivariant map algebra case, though we include further details where our arguments differ from those of \cite{NS}.  We close the section with Theorem \ref{keythmext} and Proposition \ref{singleorbit}, which give a precise description of $\Ext_\cL^1(V,W)$ for simple evaluation modules $V$ and $W$, in terms of homomorphisms of modules over finite-dimensional reductive Lie algebras.
%The results of Section 5 are entirely new.  Apart from special cases, it was not known how to calculate extensions between simple evaluation and non-evaluation modules, even for equivariant map algebras.  The block decomposition in \cite{NS} was obtained only under the key assumption that $\Ext_\cL^1(V,k_\lambda)=0$ for simple evaluation modules $V$ and 1-dimensional non-evaluation modules $k_\lambda$.  Such algebras $\cL$ are said to be {\em extension-local}.
It has remained an open question as to whether there exist equivariant map algebras which are {\em not} extension-local.  In Theorem \ref{extension-locality}, we completely settle this question and show that every Lie subalgebra of $\fg\ot_k S$ is extension-local.  In particular, twisted current algebras (and thus equivariant map algebras) are always extension-local.  It then easily follows that there are no non-trivial extensions between simple evaluation and non-evaluation modules (Theorem \ref{usingextlocal}).

%In Section 4, we return to the problem of determining extensions between arbitrary pairs of simple modules in $\LL$-$\mathbf{mod}$, and describing the corresponding block decomposition of the category.  In Theorem \ref{blockdecomp}, we obtain a bijection between blocks and pairs $(\chi,k_\alpha)$ of spectral characters $\chi$ and certain non-evaluation representations $k_\alpha$.  While the statements and proofs in this section are mostly generalizations of previous work in the special case of equivariant map algebras, we are able to offer a more complete classification by using the results of Section 5 to drop the hypothesis of extension-locality.

In Section 4, we apply our extension-locality results to the context of calculating extensions and blocks of the category $\LL$-$\mathbf{mod}$ for any twisted current algebra $\LL$.   In Theorem \ref{blockdecomp}, we obtain a bijection between blocks and pairs $(\chi,k_\alpha)$ of spectral characters $\chi$ and certain non-evaluation representations $k_\alpha$.  The results generalize previous work done in the case of equivariant map algebras \cite{NS}, and we are able to offer a more complete classification even in this case, by using the results of Section 3 to drop the hypothesis of extension-locality.

We close the paper with a specialization to the previously unexplored setting of twisted forms of current algebras.  Since these algebras are perfect, it is obvious that there are no non-evaluation modules, and the classification of blocks reduces to the calculation of spectral characters.  We show that these spectral characters may be interpreted as maps from the maximal spectrum $\Max S$ to the fundamental group $P/Q$ of the root system of $\fg$.  This generalizes previously known results for loop algebras \cite{CM,Se}.  None of this material was previously known in the general setting of twisted current algebras.

Our proofs in Section 5 are completely different from those appearing in the special case of loop algebras \cite{CM,Se}, however.  Chari, Moura, and Senesi made extensive use of Weyl modules and Drinfeld polynomials, while our approach relies on a more explicit determination of extensions, whose dimension depends on particular Littlewood-Richardson coefficients, which are then calculated with the PRV formula.  We illustrate the power of our results with a strikingly simple parametrisation of the blocks of the mysterious Margaux algebra, a twisted form of a multiloop algebra whose existence is only known cohomologically through descent theory.

\bigskip

\noindent
{\bf Acknowledgements.} The authors thank Georgia Benkart for helpful comments on a previous draft of this paper.

\bigskip

\noindent
{\bf Notation.} Throughout this paper, $k$ will denote an algebraically closed field of characteristic zero.  All algebras and tensor products will be taken over $k$ unless otherwise indicated.  The category of finitely generated unital commutative associative $k$-algebras will be denoted by {\bf $k$-alg}. We write $\Max\,S$ for the maximal spectrum of each object $S$ in {\bf $k$-alg}.  Unless otherwise specified, $L$ will denote an arbitrary Lie algebra over $k$, $L'=[L,L]$ will be its derived subalgebra, and {\bf $L$-mod} will be the category of finite-dimensional $L$-modules.  For linear forms $\lambda\in (L/L')^*$, the associated $1$-dimensional $L$-module will be denoted by $k_\lambda$, and $V^L$ will be the (trivial) submodule $V^L=\{v\in V\ :\ x.v=0\hbox{\ for all\ }x\in L\}$ of $L$-invariants of any $L$-module $V$.

%\newpage
%%%

%\ja{
%\section*{Acknowledgements}
%J.A. acknowledges ....
%M.L. would like to thank ...
%}

\section{Representations of Twisted Current Algebras}

Let $\fg$ be a finite-dimensional reductive Lie algebra over $k$, $S$ an object of {\bf $k$-alg}, and $\Gamma$ a finite group acting by $k$-algebra automorphisms on $S$.  Given a crossed homomorphism $u:\ \Gamma\rightarrow\Aut_{S-Lie}(\fg\ot S)$, $\gamma\mapsto u_\gamma$, there is an action of $\Gamma$ by $k$-Lie algebra automorphisms on the {\em current algebra} $\fg\ot S$:
$${}^{\gamma}(x\ot s)=u_\gamma(x\ot{}^{\gamma}s),$$
for all $\gamma\in\Gamma$ and $x\ot s\in \fg\ot S$.  The Lie subalgebra
$$\cL=(\fg\ot S)^\Gamma=\{z\in\fg\ot S\ |\ {}^\gamma z=z\hbox{\ for all\ }\gamma\in\Gamma\}$$
of $\Gamma$-invariants is called a {\em twisted current algebra}.  When $\fg$ is simple and $S$ is reduced, every finite group action by $k$-Lie algebra automorphisms on $\fg\ot S$ is of this form, and the twisted current algebras are precisely the invariant subalgebras of $\fg\ot S$ under finite group actions \cite[Theorem 2.2]{twcurr}.

When $\fg$ is simple and $S$ is a Galois extension of the subalgebra $S^\Gamma$ of $\Gamma$-invariants, we obtain the {\em $S/S^\Gamma$-twisted forms} of the current algebra $\fg\ot S^\Gamma$, those $S^\Gamma$-Lie algebras $\mathcal{A}$ with the property that $\mathcal{A}\ot_{S^\Gamma} S$ is isomorphic to the $S$-Lie algebra $\fg\ot S$.  See \cite{repforms} for details.  Another important class of examples occurs when the automorphisms
$u_\gamma:\ \fg\ot S\rightarrow\fg\ot S$ leave the subalgebra $\fg\ot k\cong\fg$ setwise invariant. Such twisted current algebras are called {\em equivariant map algebras}.

%Such a splitting of the action also occurs when $(\fg\ot S)^\Gamma$ is an equivariant map algebra \cite{NS}.  This is precisely the case when the image of the $1$-cocycle $u$ is in the space $(\Aut_{k-Lie}\fg)\ot 1$ of $S$-linear extensions of $k$-Lie algebra automorphisms of $\fg$.  Descent theory ensures that such a splitting also appears in the setting of twisted forms, where $\fg$ is simple and $S$ is a Galois extension of the subalgebra $S^\Gamma$ of $\Gamma$-invariants \cite{repforms}.  For the rest of this paper, we will assume the existence of a splitting of the $\Gamma$-action as in (\ref{rem:splitting}), in order to have access to the tools we need, while maintaining sufficient generality to include equivariant map algebras and twisted forms, as well as all twisted current algebras for $\fg$ simple and $S$ reduced.

The finite-dimensional simple modules of twisted current algebras $\cL=(\fg\ot S)^\Gamma$ were classified in \cite{twcurr}.  We summarize the necessary definitions and classification results in this section.

\bigskip

For each $s\in S$ and $M\in\Max\,S$, let $s(M)\in k$ be the reduction of $s$ modulo $M$:
$$s(M)+M=s+M\in S/M\cong k.$$
Likewise, we write $z(M)=\sum s_i(M)x_i\in\fg,$
for all $z=\sum x_i\ot s_i\in \fg\ot S$. For each $M\in\Max\,S$, the image of the map
\begin{align*}
\ev_M:\ &\fg\ot S\rightarrow \fg\\
&z\mapsto z(M),
\end{align*}
restricted to the twisted current algebra $\cL=(\fg\ot S)^\Gamma$, is the reductive Lie algebra
$$\fg^M=\{x\in\fg\ |\ \gamma(M)x=x\hbox{\ for all\ }\gamma\in\Gamma^M\},$$
where $\Gamma^M=\{\gamma\in\Gamma\ |\ {}^\gamma M=M\},$
and $\gamma(M)$ is the $k$-Lie algebra automorphism
\begin{align*}
\gamma(M):\ &\fg\rightarrow \fg\\
&x\mapsto ({}^\gamma(x\ot 1))(M).
\end{align*}
For any Lie subalgebra $L$ of $\fg\ot S$, the modules whose actions factor through an {\em evaluation homomorphism}
$$\ev_{\bf M}:\ L\rightarrow\ev_{M_1}(L)\oplus\cdots\oplus\ev_{M_r}(L)$$
for some ${\bf M}=\{M_1,\ldots,M_r\}\subseteq\Max\,S$ are called {\em evaluation modules}.
\begin{theorem}{\em \bf \cite[Theorems 3.1 and 3.2]{twcurr}}\label{completelist}  Let $\rho:\ \cL\rightarrow\hbox{End}_k V$ be an irreducible finite-dimensional representation of the twisted current algebra $\cL$.  Then there exist $M_1,\ldots,M_r\in\Max\,S$ in distinct $\Gamma$-orbits, together with finite-dimensional simple $\fg^{M_i}$-modules $(V_i,\rho_i)$ and a (possibly trivial) $1$-dimensional $\cL$-module $W$, such that $V\cong W\ot V_1\otimes\cdots\ot V_r$, with $\cL$-action $\rho_i\circ\ev_{M_i}$ on each $V_i$.  Conversely, every $\cL$-module of this form is finite-dimensional and simple.
\end{theorem}\qed

The classification of simple modules up to isomorphism is given in terms of a trivial fibre bundle
\begin{equation} \label{bundleiso}
\mathcal{R} \; = \bigsqcup_{M \in \Max S} \Rep(\fg^M) \twoheadrightarrow \Max S,
\end{equation}
%$$\mathcal{B}=\coprod_{M\in\Max\,S}\Rep(\fg^M)\rightarrow\Max\,S,$$
where the fibre over $M\in\Max\,S$ is the set $\Rep(\fg^M)$ of isomorphism classes of finite-dimensional simple $\fg^M$-modules.  The class $[k_0]$ of the trivial $1$-dimensional module $k_0$ is a member of each fibre.  For any section $f:\ \Max\,S\rightarrow\mathcal{R}$ of this bundle, we define the {\em support} of $f$ to be
$\supp f=\{M\in\Max\,S\ |\ f(M)\neq [k_0]\}.$
If $\supp f$ is of finite cardinality, we say that $f$ is {\em finitely supported}, and we denote the set of finitely supported sections by $\mathcal{F}$.

Given an evaluation module $\displaystyle{E=\bigotimes_{i-1}^r E_i}$, with $[(E_i,\rho_i)]\in\Rep(\fg^{M_i})$, we define a section $f_E\in\F$ as follows:
$$f_E:\ M\mapsto\left\{\begin{array}{ll}
[(E_i,\rho_i\circ(\gamma(M))^{-1})] & \hbox{if \ }M={}^\gamma M_i\hbox{\ for some\ }\gamma\in\Gamma\hbox{\ and\ }1\leq i\leq r\\
\lbrack k_0\rbrack & \hbox{otherwise,}
\end{array}\right.$$
where ${}^{\gamma}M_i$ denotes the image of the maximal ideal $M_i$ under the action of $\gamma\in\Gamma$ on $S$.  Such sections $f_E$ are invariant under an action $f\mapsto{}^\gamma f$ of the group $\Gamma$ on the set $\F$, where
\begin{equation} \label{actionisobundle}
({}^\gamma f)(M)=\left[\Big(V_{{}^{\gamma^{-1}}M},\rho_{{}^{\gamma^{-1}}M} \circ \big(\gamma(M)\big)^{-1}\Big)\right],
\end{equation}
and $(V_{{}^{\gamma^{-1}}M},\rho_{{}^{\gamma^{-1}}M})$ is a representative of the isomorphism class $f({}^{\gamma^{-1}} M)$ of $\fg^{{}^{{\gamma}^{-1}}M}$-modules. In the remainder of the paper, we will use the more concise notation $f({}^{\gamma^{-1}}M) \circ \big(\gamma(M)\big)^{-1}$ for $\big({}^\gamma f\big)(M)$.

Conversely, for each $f$ in the set $\F^\Gamma$ of finitely supported $\Gamma$-invariant sections, we can decompose the support of $f$ into pairwise disjoint $\Gamma$-orbits:
$$\supp f=\bigsqcup_{i=1}^r \Gamma.M_i,$$
and define an $\cL$-module:
$$V_f=V_{f(M_1)}\ot\cdots\ot V_{f(M_r)},$$
where the $V_{f(M_i)}$ are representatives of the isomorphism classes $f(M_i)\in\Rep(\fg^{M_i})$, and $(V_f,\rho_f)$ is an evaluation $\cL$-module via the pullback of $\ev_{\bf M}$.  The isomorphism class of $(V_f,\rho_f)$ is independent of the choice of representatives $V_{f(M_i)}$ and these give a natural bijection between $\F^\Gamma$ and the set $\mathcal{S}$ of simple objects in $\cL_{\ev}-\mathbf{mod}$, the full subcategory of finite-dimensional evaluation modules in $\cL-\mathbf{mod}$.
\begin{theorem}{\bf \cite[Proposition 3.13]{twcurr}} The map $f\mapsto V_f$ is a bijection from $\F^\Gamma$ to $\mathcal{S}$, with inverse $E\mapsto f_E$.
\end{theorem}

Using this equivalence, we define the {\em support} of a simple object $E\in\cL_{\ev}-\mathbf{mod}$ to be the support of the corresponding section $f_E$.  The support of an arbitrary module in $\cL_{\ev}-\mathbf{mod}$ is defined as the union of the supports of its simple subquotients.  A finite-dimensional simple evaluation module supported on a single $\Gamma$-orbit $\Gamma.M$ is called a {\em single-orbit evaluation module}, and the full subcategory of such modules in $\cL_{\ev}-\mathbf{mod}$ is denoted by $\cL_{\ev}^{\Gamma.M}-\mathbf{mod}$. By \cite[Lemma 3.8]{twcurr}, there are natural equivalences of categories
\begin{eqnarray} \label{catequiv}
\LL_\ev^{\Gamma.M}\text{-}\mathbf{mod} & \cong & \LL_\ev^{\Gamma.({}^\gamma M)}\text{-}\mathbf{mod } =  \LL_\ev^{\Gamma.M}\text{-}\mathbf{mod} \\
(\rho,V) & \longmapsto & \left(\rho \circ \big(\gamma(M)\big)^{-1},V\right), \notag
\end{eqnarray}
for all $M\in\Max\,S$ and $\gamma\in\Gamma$.

The classification of finite-dimensional simple $\cL$-modules up to isomorphism finishes with a removal of redundancies created by the tensor product with the $1$-dimensional module $W$ in the statement of Theorem \ref{completelist}:

\begin{theorem}{\bf \cite[Theorem 3.14]{twcurr}}\label{threeconditions} The isomorphism classes of finite-dimensional simple $\cL$-modules are in natural bijection with the pairs $(\lambda,f)\in\cL^*\times\F^\Gamma$, such that $\lambda$ vanishes on $[\cL,\cL]$, $\cL/\ker{\rho_f}$ is semisimple, and $\ker(\lambda\ot 1+1\ot\rho_f)=(\ker\lambda)\cap(\ker\rho_f)$.  Explicitly, $(\lambda,f)$ corresponds to the isomorphism class of $k_\lambda \ot V_f$.
\end{theorem}

\section{Extensions between Evaluation and Non-Evaluation \\Modules}

Given an abelian category $\mathcal{C}$, two indecomposable objects $I_1$ and $I_2$ are {\em linked}, and we write $I_1\sim I_2$ if there is no decomposition $\mathcal{C} = \mathcal{C}_1 \oplus \mathcal{C}_2$ where $\mathcal{C}_i$ are abelian subcategories such that $I_i \in \mathcal{C}_i$. For a given linkage class $\Lambda$, define $\mathcal{C}_\Lambda$ to be the full subcategory of $\mathcal{C}$ consisting of direct sums of indecomposables from $\Lambda$. The category $\mathcal{C}_\Lambda$ is called a {\em block} of $\mathcal{C}$ and the category $\mathcal{C}$ then has the following {\em block decomposition}:  $\mathcal{C} = \bigoplus \mathcal{C}_\Lambda$, where the direct sum is taken over all linkage classes $\Lambda$.

The blocks of the category of finite-dimensional representations of equivariant map algebras were classified in \cite{NS}, modulo the hypothesis of extension-locality, defined as below.

\begin{definition}{\em
A Lie subalgebra $\cL$ of the current algebra $\fg\ot S$ is {\em extension-local} if $\Ext_\cL^1(k_\lambda,V)=0$ whenever $V$ is a simple evaluation module and $k_\lambda$ is not an evaluation module.\footnote{The original definition was that $\Ext_\cL^1(V,k_\lambda)=0$ for all simple $V\in\LL_\ev$-$\mathbf{mod}$ and $k_\lambda\notin\LL_\ev$-$\mathbf{mod}$.  Since $\hbox{Ext}_\cL^1(V^*,k_\lambda^*)=\hbox{H}^1(\cL;V\ot k_\lambda^*)=\hbox{Ext}_\cL^1(k_\lambda,V)$, this is clearly equivalent to the definition above, since duals of (non)evaluation modules are (non)evaluation.}}
\end{definition}

Despite a number of interesting examples of extension-local equivariant map algebras, it has remained an open question whether there exist equivariant map algebras which are {\em not} extension-local \cite[Remark 5.14]{NS}.  In this section, we prove that there are no such examples; every twisted current algebra is extension-local, and the hypothesis of extension-locality is thus superfluous in \cite{NS} and in the present article.

\begin{theorem}\label{extension-locality}
Let $\cL$ be a Lie subalgebra of $\fg\ot S$.  Then $\cL$ is extension-local.  In particular, every twisted current algebra is extension-local.
\end{theorem}

\noindent
\proof  Suppose $\cL$ is not extension-local.  Then there is a simple evaluation representation $(V,\psi)$ and a 1-dimensional non-evaluation $\cL$-module $k_\lambda$ such that $\Ext_\cL^1(k_\lambda,V)\neq 0$.  As usual, we view $\lambda$ as a linear functional on $\cL$, with kernel containing the derived subalgebra $\cL'$.  In particular, there exist maximal ideals $M_1,\ldots,M_r\in\Max\,S$ such that
$$\ker\ev_{\mathbf{M}}:=\bigcap_{i=1}^r\ker\ev_{M_i}$$
lies in $\ker\psi$.  As $k_\lambda$ is not an evaluation module, $\ker\ev_{\mathbf{M}}\not\subseteq\ker\lambda$, so there exists
\begin{equation}\label{4etoiles}
z\in\ker\ev_{\mathbf{M}}\setminus\ker\lambda.
\end{equation}
Rescaling $z$ if necessary, we may assume that $\lambda(z)=-1$.

As $\Ext_\cL^1(k_\lambda,V)\neq 0$, there is an $\cL$-module structure
$$\phi:\ \cL\rightarrow\End_k(V\oplus k_\lambda)$$
on the vector space direct sum $V\oplus k_\lambda$, such that the natural inclusion and projection maps give a nonsplit short exact sequence
\begin{equation}\label{1etoile}
0\rightarrow V\rightarrow V\oplus k_\lambda\rightarrow k_\lambda\rightarrow 0.
\end{equation}
Fixing a $k$-basis $\{v_1,\ldots,v_\ell\}$ for the $\cL$-module $V$ and a nonzero element $v'\in k_\lambda$, we may express the image $\phi(x)$ of each $x\in\cL$ as an $(\ell+1)\times(\ell+1)$ matrix
$$\phi(x)=\left(\begin{array}{cc}
A(x) & \bb(x)\\
0 & \lambda(x)
\end{array}\right),$$
with respect to the basis $\{v_1,\ldots,v_\ell,v'\}$.  Here $A(x)$ is an $\ell\times\ell$ matrix in $M_\ell(k)$ and $\bb(x)$ is an $\ell\times 1$ column vector in $k^\ell$.

\bigskip

For clarity, we will divide the rest of the proof into a series of short steps.

\bigskip

\noindent
{\it Step 1.} Suppose that the $k$-vector space
$$\mathcal{R}=\left\{\left(\begin{array}{c}
\bb(x)\\
\lambda(x)
\end{array}\right)\ |\ x\in\cL\right\}$$
is 1-dimensional.  Then $V$ is 1-dimensional.

\medskip

\noindent
{\it Proof.}  By hypothesis, there exists $\bb\in k^\ell$ such that
$$\bb(x)=\lambda(x)\bb$$
for all $x\in\cL$.  If $x\in\cL$ and $y\in\ker\lambda$, we have
$$\phi[x,y]=[\phi(x),\phi(y)]=\left(\begin{array}{cc}
[A(x),A(y)]& -\lambda(x)A(y)\bb\\
0& 0
\end{array}
\right).$$
In particular,
$$-\lambda(x)A(y)\bb=\bb([x,y])=\lambda([x,y])\bb=0,$$
since $\lambda$ vanishes on $\cL'$.  Specializing to any $x\not\in\ker\lambda$, we see that
\begin{equation}\label{2etoiles}
A(y)\bb=0
\end{equation}
for all $y\in\ker\lambda$.

If $\bb=0$, then
$$\phi(x)=\left(\begin{array}{cc}
A(x)&0\\
0&\lambda(x)
\end{array}\right)$$
for all $x\in\cL$, which is impossible since (\ref{1etoile}) is not split.

We can thus identify $\bb$ with a nonzero element $\sum_{i=1}^\ell b_iv_i$ of $V$.  For $x\in\cL$ and $y\in\ker\lambda$,
\begin{equation}\label{3etoiles}
y.(x.\bb)=[y,x].\bb+x.(y.\bb).
\end{equation}
Since $y,[y,x]\in\ker\lambda$ and the action of each element $u\in\cL$ on $\bb$ is given by $u.\bb=A(u)\bb,$ we see that $y.(x.\bb)=0$,
by (\ref{2etoiles}) and (\ref{3etoiles}).  Moreover, $\mathcal{U}(\cL).\bb=V$ by the irreducibility of $V$, so
$$0=(\ker\lambda).(\mathcal{U}(\cL).\bb)=(\ker\lambda).V.$$
Thus $\ker\lambda\subseteq\ker\psi$, so $\cL/\ker\psi$ is of dimension at most $1$, and the irreducible $\cL$-module $V$ is 1-dimensional.

\bigskip

\noindent
{\it Step 2.}  There exists $t\in\ker\lambda$ such that $\bb(t)\neq 0$.

\medskip

\noindent
{\it Proof.}  Suppose that $\mathcal{R}$ is 1-dimensional.  By Step 1, $V=k_\mu$ for some $\mu\in\cL^*$, and
$$\phi(x)=\left(\begin{array}{cc}
\mu(x)&\lambda(x)\bb\\
0&\lambda(x)
\end{array}\right)$$
for some nonzero $\bb\in k$.  Moreover,
$$\mu(y)\bb=A(y)\bb=0$$
for all $y\in\ker\lambda$ by (\ref{2etoiles}), so $\ker\lambda\subseteq\ker\mu$ and $\mu=\lambda$ or $\mu=0$.  As $V=k_\mu$ is an evaluation module and $k_\lambda$ is not, we see that $\mu=0$.

Therefore,
$$\phi(x)=\left(\begin{array}{cc}
0&\lambda(x)\bb\\
0&\lambda(x)
\end{array}
\right)$$
for all $x\in \cL$.  But then the map
\begin{eqnarray}
\rho:\ k_\lambda &\longrightarrow&V\oplus k_\lambda\label{2csillagok}\\
cv'&\longmapsto&(c\bb,cv')\nonumber
\end{eqnarray}
for each $cv'\in k_\lambda=\hbox{Span}\{v'\}$ is a splitting of (\ref{1etoile}), contradicting the assumption that the exact sequence (\ref{1etoile}) is not split.  Hence $\mathcal{R}$ is at least $2$-dimensional.

We can thus construct $t\in\ker\lambda$ such that $\bb(t)\neq 0$, by taking an appropriate linear combination of any pair $u_1,u_2\in\cL$ such that
$$\left(\begin{array}{c}
\bb(u_1)\\
\lambda(u_1)
\end{array}\right)\hbox{\quad and\quad}\left(\begin{array}{c}
\bb(u_2)\\
\lambda(u_2)
\end{array}\right)$$
are linearly independent.

\bigskip

\noindent
{\it Step 3.}  There exists $w\in\ker\ev_{\mathbf{M}}\cap\ker\lambda$ such that $\phi(w)\neq 0$.

\medskip

\noindent
{\it Proof.}  Let $X=\{\phi(x)\ :\ x\in\ker\ev_{\mathbf{M}}\}$.  If $X=0$, then $\ker\ev_{\mathbf{M}}\subseteq\ker\phi\subseteq\ker\lambda$.  Since $k_\lambda$ is not an evaluation module, this is impossible.  Hence $X\neq 0$, and without loss of generality, we can assume that each nonzero $\phi(x)\in X$ is of the form
$$\phi(x)=\left(\begin{array}{ll}
0&\bb(x)\\
0&\lambda(x)
\end{array}\right),$$
with $\lambda(x)\neq 0$.  In particular, $X$ is $1$-dimensional, as
$$\phi(\lambda(y)x-\lambda(x)y)=\left(\begin{array}{cc}
0&\lambda(y)\bb(x)-\lambda(x)\bb(y)\\
0&0
\end{array}\right)$$
will be nonzero for any pair of linearly independent $\phi(x),\phi(y)\in X$.  There is thus a nonzero element $\bb\in V=k^\ell$ such that
$$\phi(x)=\left(\begin{array}{cc}
0&\lambda(x)\bb\\
0&\lambda(x)
\end{array}\right),$$
for all $x\in\ker\ev_{\mathbf{M}}$.

For each $x\in\ker\ev_{\mathbf{M}}$ and $y\in\ker\lambda$, we have
$$\phi([x,y])=[\phi(x),\phi(y)]=\left(\begin{array}{cc}
0&-A(y)\bb(x)-\lambda(x)\bb(y)\\
0&0
\end{array}\right).$$
If $-A(y)\bb(x)-\lambda(x)\bb(y)$ is nonzero for some $x\in\ker\ev_{\mathbf{M}}$ and $y\in\ker\lambda$, then we are done, since $w=[x,y]\in\ker\ev_{\mathbf{M}}\cap\ker\lambda$ and $\phi(w)\neq 0$.

We can thus assume that
\begin{equation}\label{1csillag}
A(y)\lambda(x)\bb+\lambda(x)\bb(y)=A(y)\bb(x)+\lambda(x)\bb(y)=0,
\end{equation}
for all $x\in\ker\ev_{\mathbf{M}}$ and $y\in\ker\lambda$.  The map $\rho:\ k_\lambda\rightarrow V\oplus k_\lambda$ in (\ref{2csillagok}) is then again a splitting of the short exact sequence (\ref{1etoile}).  Indeed, since $\ker\lambda$ is codimension $1$ in $\cL$ and $\ker\ev_{\mathbf{M}}\not\subseteq\ker\lambda$, we see that $\ker\ev_{\mathbf{M}}+\ker\lambda=\cL$, and we need only verify that $\rho$ commutes with the actions of $\ker\ev_{\mathbf{M}}$ and $\ker\lambda$.

Let $cv'\in k_\lambda=\hbox{Span}\{v'\}$.  For $x\in\ker\ev_{\mathbf{M}}$, we have
\begin{equation}\label{3csillagok}
x.\rho(cv')=x.(c\bb,cv')=(\lambda(x)c\bb,\lambda(x)cv').
\end{equation}
Likewise, for $x\in\ker\lambda$, we have
\begin{equation}\label{4csillagok}
x.\rho(cv')=x.(c\bb,cv')=(cA(x)\bb+c\bb(x),0).
\end{equation}
Taking $z$ as in (\ref{4etoiles}), we see that $\lambda(-cz)=c$, so
$x.\rho(cv')=(0,0)=(\lambda(x)c\bb,\lambda(x)cv')$ by (\ref{1csillag}) for all $x\in\ker\lambda$.  In both (\ref{3csillagok}) and (\ref{4csillagok}),
$$\rho(x.cv')=\rho(\lambda(x)cv')=(\lambda(x)c\bb,\lambda(x)cv')=x.\rho(cv'),$$
so the map $\rho$ is a splitting of (\ref{1etoile}).  This contradicts the non-triviality of the extension, so there exists $w\in\ker\ev_{\mathbf{M}}\cap\ker\lambda$ with $\phi(w)\neq 0$.

\bigskip

\noindent
{\it Step 4.}  $\Ext_\cL^1(k_\lambda,V)=0$.

\medskip

\noindent
{\it Proof.}  By Step 3, there exists $w\in\cL$ such that
$$\phi(w)=\left(\begin{array}{cc}
0&\bb(w)\\
0&0
\end{array}\right)\neq 0.$$
Choosing $z$ as in (\ref{4etoiles}), we have
$$\phi([z,w])=[\phi(z),\phi(w)]=\left[\left(\begin{array}{cc}
0&\bb(z)\\
0&-1
\end{array}\right),\left(\begin{array}{cc}
0&\bb(w)\\
0&0
\end{array}\right)\right]\\
=\phi(w).$$
Let $z_n=(\hbox{ad}\,z)^{n-1}w$ for all $n\geq 1$.  Then
$$\phi(z_n)=\phi(w)\neq 0$$ for all $n$.  As $z,w\in\ker\ev_{\mathbf{M}}\subseteq\cL\cap(\fg\ot M_1)$, we also have
$$z_n\in\cL\cap(\fg\ot M_1^n),$$
where $M_1^n$ is the $n$th power of the maximal ideal $M_1\in\hbox{Max}\,S$.

Let $\mathcal{I}_n=(\cL\cap(\fg\ot M_1^n))+\ker\phi$.  Clearly,
$$\dim_k(\mathcal{I}_n/\ker\phi)\leq \dim_k(\cL/\ker\phi)\leq (\ell+1)^2.$$
That is,
$$\mathcal{I}_1/\ker\phi\ \supseteq\  \mathcal{I}_2/\ker\phi\ \supseteq\ \cdots$$
is a decreasing sequence of finite-dimensional vector spaces.  In particular, there exists $N>0$ such that
$$\mathcal{I}_N=\mathcal{I}_{N+1}=\mathcal{I}_{N+2}=\cdots.$$
Since $z_N\in\mathcal{I}_N\setminus\ker\phi$, we see that $\dim_k(\mathcal{I}_N/\ker\phi)\neq 0,$ so
$$\left(\cL\cap\left(\fg\ot\bigcap_{n>0}M_1^n\right)\right)+\ker\phi=\bigcap_{n>0}\mathcal{I}_n=\mathcal{I}_N$$
is {\it not} contained in $\ker\phi$.

The finitely generated $k$-algebra $S$ is a quotient of a polynomial ring $T$ in finitely many variables.  The maximal ideal $M_1\in\Max\,S$ is thus the quotient of a maximal ideal $N$ of $T$.  As $T$ is a noetherian domain,
$$\bigcap_{n>0}N^n=0,$$
by the Krull Intersection Theorem.  It follows that
$$\bigcap_{n>0}M_1^n=0\hbox{\quad and\quad}\bigcap_{n>0}\mathcal{I}_n=\ker\phi,$$
a contradiction.  Hence $\Ext_\cL^1(k_\lambda,V)=0$.\qed

%%%
%\include{newsectionIV.tex}
%%%
\section{Blocks for Twisted Current Algebras}

% NEW PLAN
%%%%%%%%%
% I. Intro
% II. Background
% III. Ext-locality
% IV. Application to twisted current algebras
% V. Specialization to twisted forms

% PLAN FOR THIS NEW SECTION...
%%%%%%%%%%%%%%%%%%%%%%%%%
% Basically, we apply the ext-locality result to make some progress.
%
% Contents:
% --> Main results of extensions between modules (including non-evaluation ones)
% --> Main results of block decompositions

In this section, we compute extensions in the category $\LL$-$\mathbf{mod}$ of finite-dimensional modules over a twisted current algebra $\LL = (\fg \otimes S)^\Gamma$.  We use Theorem \ref{extension-locality} to reduce the calculation of ext-blocks to the calculation of extensions between simple evaluation modules, and then use the description of ext-blocks to determine the blocks for $\LL$-$\mathbf{mod}$.  The material in this section relies on basic definitions and facts about extensions found in Appendix A, and extends previously known results for equivariant map algebras \cite{NS} to the general setting of twisted current algebras.

%We combine our results on extension-locality (Section 3) with the classification of simple objects \cite{twcurr} and generalizations of the equivariant map algebra case \cite{NS} to describe the blocks of  $\LL$-$\mathbf{mod}$. %The case of equivariant map algebras, as presented in \cite{NS}, and the description of simple finite-dimensional $\LL$-modules \cite{twcurr} are the two main ingredients for what follows.

%To obtain our description,  we use extension-locality to reduce the problem to a calculation of extension blocks of $\LL_\ev$-$\mathbf{mod}$. We then describe the sets $\Ext_\LL^1(E,F)$ for simple objects $E$ and $F$ of the category $\LL_\ev$-$\mathbf{mod}$ of finite-dimensional evaluation representations, and use this information to classify extension blocks. We conclude the section with a classification of the blocks of both $\LL_\ev$-$\mathbf{mod}$ and $\LL$-$\mathbf{mod}$.

We recall the definition of extension block and its relation to block decomposition.
%the sets $\Ext_\LL^1(V,W)$ for simple modules $V$ and $W$ in $\LL$-$\mathbf{mod}$. We then proceed to describing the extension blocks of $\LL$-$\mathbf{mod}$.
%In this section, we describe the sets $\Ext_\LL^1(V,W)$ and for finite-dimensional modules $V$ and $W$ over a twisted current algebra $\LL = (\fg \otimes S)^\Gamma$.
%%The main ingredients are arguments of Neher and Savage in the equivariant map setting \cite{NS}, and the classification of simple finite-dimensional $\LL$-modules by one of the authors of the present paper \cite{twcurr}.
%%Our work is inspired by analogous results of Neher and Savage \cite{NS} in the special case where $\LL$ is an equivariant map algebra, and the classification of simple finite-dimensional $\LL$-modules by one of the authors of the present paper \cite{twcurr}.
%The case of equivariant map algebras, as presented in \cite{NS}, and the description of simple finite-dimensional $\LL$-modules \cite{twcurr} are the two main ingredients for what follows.

%Let $\LL$ be a twisted current algebra. %and consider the categories $\LL$-$\mathbf{mod}$ and $\LL_\ev$-$\mathbf{mod}$.  % of finite-dimensional $\LL$-modules and finite-dimensional $\LL$ evaluation modules, respectively.
%The blocks of $\LL$-$\mathbf{mod}$ and $\LL_\ev$-$\mathbf{mod}$ may be understood in terms of extensions of the simple objects they contain.
%In this section, we first relate extensions between simple objects to blocks, and then we describe the blocks of the categories $\LL$-$\mathbf{mod}$ and $\LL_\ev$-$\mathbf{mod}$.
% % % >> GENERAL BLOCKS << % % %
\begin{definition} \label{defextblock}\label{blockdef}{\em Let $V$ and $W$ be simple objects of an abelian category $\mathcal{C}$ of finite length, such as $\LL$-$\mathbf{mod}$, $\LL_\ev$-$\mathbf{mod}$, or the finite-dimensional modules over a finite-dimensional Lie algebra.  We say that $V$ and $W$ are in the same {\em extension block} of $\mathcal{C}$ if $V \cong W$ or there is a finite sequence of simple objects $\{T^i\}_{i=0}^N$ with $T^0 = V$, $T^N = W$, and
$$ \Ext_\mathcal{C}^1(T^i,T^{i+1}) \neq 0 \qquad \qquad \text{or} \qquad \qquad  \Ext_\mathcal{C}^1(T^{i+1},T^i) \neq 0 $$
for every $i \in \{0,\,..., N-1\}$. The set of extension blocks of $\mathcal{C}$ is denoted by $\Bloc(\mathcal{C})$, and is a partition of the collection of simple objects of $\mathcal{C}$ into disjoint subcollections.  The block containing a simple object $V$ will be denoted by $\llbracket V \rrbracket$.

Let $b$ be an extension block of $\mathcal{C}$.  The full subcategory $\mathcal{C}^b$ of objects whose simple subquotients are all contained in $b$ is a {\em block} of the category $\mathcal{C}$, as defined at the beginning of Section 3.  Conversely, every block of $\mathcal{C}$ is of this form.  More details may be found in \cite[Section 1.13]{Humphreys-CatO}, for instance.}
%An object of $\mathcal{C}$ is said to be in the {\em block} $\mathcal{C}^b$ determined by $b$ if each of its simple subquotients is in $b$.  Each $\mathcal{C}^b$ is a full subcategory of $\mathcal{C}$, with no non-trivial extensions between objects in distinct $\mathcal{C}^b$.  It is clear from the definition that the blocks of the category $\mathcal{C}$ are in bijection with the extension blocks of $\mathcal{C}$.}
\end{definition}

%\begin{remark} {\em When the category is clear from the context, the extension block of a simple object $V$ will be denoted by $\llbracket V \rrbracket$.} %It is worth noting that $\, W \cong V \, \Rightarrow \, \llbracket W \rrbracket = \llbracket V \rrbracket$.
%\end{remark}

%\begin{proposition} \label{blockdef} Let $\mathcal{C}$ be as in Definition \ref{defextblock}. Then the blocks of $\mathcal{C}$ are precisely the following:
%$$ \mathcal{C}^b \;\; = \;\; \left\{\;\; V \in {\it ob}(\mathcal{C}) \quad \Bigg| \; \begin{array}{c} \text{All simple subquotients of } V \\
%\text{are in the extension block } b. \end{array}\right\} $$
%where $b \in \Bloc(\mathcal{C})$. In particular, the blocks of the category $\mathcal{C}$ are in bijection with its extension blocks $\Bloc(\mathcal{C})$.
%\end{proposition}
%\proof See \cite[\S 1.13]{Humphreys-CatO} for instance. \qed
% % %

%Proposition \ref{blockdef} shows that it is sufficient to understand the extension blocks of \linebreak $\LL$-$\mathbf{mod}$ in order to understand its blocks. Thus, it is important to identify for which combinations of simple objects $V$ and $W$ of $\LL$-$\mathbf{mod}$ do the set $\Ext_\LL^1(V,W)$ vanish or not. Since $\LL$ is always extension-local by Theorem \ref{extension-locality}, studying the extension blocks of $\LL$-$\mathbf{mod}$ is equivalent to studying those of $\LL_\ev$-$\mathbf{mod}$. The next two results explain this reasoning in more detail.
We begin our study of block decomposition with the following corollary to Theorem \ref{extension-locality}, which says that there are no non-trivial extensions between simple evaluation and non-evaluation modules.

\begin{theorem} \label{usingextlocal} Let $k_\lambda \otimes V_f$ and $k_\mu \otimes V_g$ be simple objects of $\LL$-$\mathbf{mod}$ where $(\lambda,f)$, $(\mu,g) \in \LL^* \times \secr^\Gamma$ satisfy the three conditions of Theorem \ref{threeconditions}. Assume that $k_{\lambda - \mu}$ is not an evaluation module. Then $\Ext_\LL^1(k_\lambda \otimes V_f,k_\mu \otimes V_g) = 0$.	
\end{theorem}
\proof Consider the following natural isomorphism
$$ \Ext^1_\LL(k_\lambda \otimes V_f,k_\mu \otimes V_g) \cong \Ext^1_\LL(k_{\lambda - \mu},V_f^*\otimes V_g) . $$
By \cite[Theorem 4.12.1]{Etingof}, the tensor product $V_f^* \otimes V_g$ is a semisimple $\LL$-module. As both $V_f$ and $V_g$ are finite dimensional evaluation modules, $V_f^* \otimes V_g$ is a finite direct sum of simple finite dimensional evaluation $\LL$-modules. Finally, since the bifunctor $\Ext_\LL^1(-,-)$ commutes with finite direct sum in both slots, the result now follows as a corollary to Theorem \ref{extension-locality}.\qed
\ \\

%Having a way to decide whether a 1-dimensional module is of evaluation type or not will be useful.

\begin{definition} {\em Define $(\LL/\LL')_\ev^* = \{ \lambda \in (\LL/\LL')^* \; | \; k_\lambda \text{ is an evaluation module}\}$. Let $(\LL/\LL')_{\nev}^*$ be a complementary vector subspace for $(\LL/\LL')_\ev^*$ in $\subseteq (\LL/\LL')^* $, so that
\begin{equation} \label{choicecomplement}
 (\LL/\LL')^* = (\LL/\LL')_{\ev}^* \oplus (\LL/\LL')_{\nev}^* \; .
\end{equation}}
\end{definition}
%We deduce the following:
\begin{corollary} \label{reducetoev} Let $k_\lambda \otimes V_f$ and $k_\mu \otimes V_g$ be simple objects of $\LL$-$\mathbf{mod}$ where $(\lambda,f)$, $(\mu,g) \in \LL^* \times \secr^\Gamma$ satisfy the three conditions of Theorem \ref{threeconditions}. Decompose $\lambda = \lambda_\ev + \lambda_{\nev}$ and $\mu = \mu_\ev + \mu_{\nev}$ according to~\eqref{choicecomplement}. Then $\llbracket k_\lambda \otimes V_f \rrbracket = \llbracket k_\mu \otimes V_g \rrbracket$ in $\LL$-$\mathbf{mod}$ % \in \Bloc(\LL\text{-}\mathbf{mod})$ %if and only if
\begin{equation*}
\qquad \Longleftrightarrow \qquad \left[ \; \begin{array}{l} \lambda_{\nev} = \mu_{\nev} \\ \llbracket k_{\lambda_\ev} \otimes V_f \rrbracket = \llbracket k_{\mu_\ev} \otimes V_g \rrbracket \text{ in } \LL_\ev\text{-}\mathbf{mod} \end{array} \right.
\end{equation*}
\end{corollary}
\proof Since $k_{\lambda - \mu}$ is not an evaluation module if an only if $\lambda_{\nev} \neq \mu_{\nev}$, Theorem \ref{usingextlocal} implies that $\lambda_{\nev} = \mu_{\nev}$ is necessary. Now suppose that $\lambda_{\nev} = \mu_{\nev}$ and consider proving the `only if' part. Without loss of generality, the problem reduces to the case where $\Ext_\LL^1(k_\lambda \otimes V_f,k_\mu \otimes V_g) \neq 0$ and the result follows from the natural isomorphism $\Ext_\LL^1(k_\lambda \otimes V_f,k_\mu \otimes V_g) \cong \Ext^1_\LL(k_{\lambda_\ev} \otimes V_f , k_{\mu_\ev} \otimes V_g)$. The `if' part follows similarly. \qed

\bigskip

% STUDY EXT BETWEEN SIMPLE EVALUATION MODULES
Next, we proceed to describing the sets $\Ext_\LL^1(E,F)$ of extensions between simple finite-dimensional evaluation $\LL$-modules $E$ and $F$ (simple objects of $\LL_\ev$-$\mathbf{mod}$). %as this will be sufficient to obtain a description of the blocks of $\LL$-$\mathbf{mod}$.

Recall that the $\LL$-action on a simple object of $\LL_\ev$-$\mathbf{mod}$ factors through the action of a reductive Lie algebra $\fg^\mathbf{M} = \bigoplus_{i=1}^r \fg^{M_i} = \ev_\mathbf{M}(\LL)$ where $\mathbf{M} = \{M_i\}_{i=1}^r$ is a set of maximal ideals of $S$ coming from $r$ distinct $\Gamma$-orbits. The reductive Lie algebra $\fg^{\mathbf{M}}$ can be neatly described in terms of $\LL$ itself and two ideals of $\LL$, $\mathfrak{Z}_{\mathbf{M}} = \ev_{\mathbf{M}}^{-1}\big(Z(\fg^\mathbf{M})\big)$ and $K_\mathbf{M} = \ker (\ev_\mathbf{M})$. The Chinese Remainder Theorem then gives an isomorphism of $S^\Gamma$-Lie algebras
\begin{equation} \label{reductivedecomp}
\fg^\mathbf{M} \cong \LL/K_\mathbf{M} \; \longrightarrow \; \LL/(\LL' + K_\mathbf{M}) \oplus \LL/\mathfrak{Z}_\mathbf{M}.
\end{equation}
In particular, we see that $\LL/\mathfrak{Z}_\mathbf{M} \cong [\fg^{\mathbf{M}},\fg^{\mathbf{M}}]$ is the derived subalgebra of the Lie algebra $\LL/(\LL' + K_\mathbf{M}) \oplus \LL/\mathfrak{Z}_\mathbf{M}$, and $\LL/(\LL' + K_\mathbf{M}) \cong Z(\fg^\mathbf{M})$ is its centre.

The following proposition and theorem can be proved with the obvious generalizations of proofs previously used in the equivariant map algebra case \cite[Propositions 3.3 and 3.6]{NS}.

%The ideals of $\LL$ appearing in \eqref{reductivedecomp} will be important for describing the extension blocks of $\LL$-$\mathbf{mod}$:

%\begin{remark} {\em If instead of one maximal ideal of $S$, we consider a set $\mathbf{M} = \{M_i\}_{i=1}^r$ of maximal ideals of $S$ coming from $r$ distinct $\Gamma$-orbits, the evaluation map
%$$ \ev_\mathbf{M} = \oplus_{i=1}^r \ev_{M_i} : \; \LL \rightarrow \bigoplus_{i=1}^r \fg^{M_i} = \fg^\mathbf{M} $$
%remains surjective, and the obvious analogue of (\ref{reductivedecomp}) holds for $\fg^\mathbf{M}$. %%Its still holds because the map $\ev_\mathbf{M}$ is (just) the ''direct sum'' of the evaluation maps $\ev_{M_i}$.
%
%%From a theoretical point of view, it might eventually help to get a (relevant) characterisation of 1-dimensional evaluation representations of $\LL$.
%}
%\end{remark}

\begin{proposition} \label{structureKab} Let $\mathbf{M} = \{M_i\}_{i=1}^r$ be a set of maximal ideals of $S$ coming from distinct $\Gamma$-orbits. Let $I_\mathbf{M} = \bigcap_{\gamma, \, i}^{\,} {}^\gamma M_i$, $I_\mathbf{M}^\Gamma = I_\mathbf{M} \cap S^\Gamma$, and
\begin{align*}
N_\mathbf{M} = \{x \in K_\mathbf{M} \, | \, I_\mathbf{M}^\Gamma.\,x \subseteq K_\mathbf{M}'\} && \tilde{N}_{\mathbf{M}} = \{y \in \mathfrak{Z}_\mathbf{M} \, | \, I_\mathbf{M}^\Gamma.\,y \subseteq \mathfrak{Z}_\mathbf{M}'\}.
\end{align*}
Then $K_\mathbf{M}' \unlhd N_\mathbf{M} \unlhd K_\mathbf{M}$ and $\mathfrak{Z}'_\mathbf{M}\unlhd \tilde{N}_\mathbf{M} \unlhd \mathfrak{Z}_\mathbf{M}$ are all $S^\Gamma$-ideals of $\LL$, with
$$ N_\mathbf{M}/K_\mathbf{M}'=\left({K_\mathbf{M}}/{K_\mathbf{M}'}\right)^{I_\mathbf{M}^\Gamma} \quad \text{and } \quad \tilde{N}_{\mathbf{M}}/\mathfrak{Z}_\mathbf{M}'=\left(\mathfrak{Z}_\mathbf{M}/\mathfrak{Z}_{\mathbf{M}}'\right)^{I_{\mathbf{M}}^\Gamma}.$$
The adjoint action of $\LL$ induces actions of $\fg^\mathbf{M} \cong \LL / K_\mathbf{M}$ on the quotients $N_\mathbf{M}/K_\mathbf{M}'$, $K_\mathbf{M}/K_\mathbf{M}'$, and $K_\mathbf{M}/N_\mathbf{M}$. Moreover,
\begin{itemize}
\item[\em(1)] $K_\mathbf{M}/N_\mathbf{M}$ is a trivial $\fg^\mathbf{M}$-module.
\item[\em(2)] As $\fg^\mathbf{M}$-modules, $N_\mathbf{M}/K_\mathbf{M}' = \bigoplus_{i=1}^r T_i$, where the $T_i$ are finite-dimensional over $k$, and $\fg^{M_j}.T_i = 0$ for all $j \neq i$. %In particular, $N_\mathbf{M}/K_\mathbf{M}'$ is finite-dimensional.
\end{itemize}
Similarly, the adjoint action of $\LL$ induces an action of $[\fg^\mathbf{M},\fg^\mathbf{M}] \cong \LL/\mathfrak{Z}_\mathbf{M}$ on the quotients $\tilde{N}_\mathbf{M}/\mathfrak{Z}_\mathbf{M}'$, $\mathfrak{Z}_\mathbf{M}/\mathfrak{Z}_\mathbf{M}'$, and $\mathfrak{Z}_\mathbf{M}/\tilde{N}_\mathbf{M}$. Moreover,
\begin{itemize}
\item[\em(3)] $\mathfrak{Z}_\mathbf{M}/\tilde{N}_\mathbf{M}$ is a trivial $[\fg^\mathbf{M},\fg^\mathbf{M}]$-module.
\item[\em(4)] As $[\fg^\mathbf{M},\fg^\mathbf{M}]$-modules, $\tilde{N}_\mathbf{M}/\mathfrak{Z}_\mathbf{M}' = \bigoplus_{i=1}^r U_i$, where the $U_i$ are finite-dimensional over $k$, and $[\fg^{M_j},\fg^{M_j}].U_i = 0$ for all $j \neq i$. %In particular, $\tilde{N}_\mathbf{M}/\mathfrak{Z}_\mathbf{M}'$ is finite-dimensional.
\end{itemize}
\end{proposition}\qed

\begin{theorem} \label{keythmext} Let $E$ and $F$ be simple objects of $\LL_\ev$-$\mathbf{mod}$ and let $\{M_i\}_{i=1}^r \subseteq \Max S$ be a set of representatives of each $\Gamma$-orbit in $\,\supp E \, \cup \, \supp F$. Write $E = \bigotimes_{i=1}^r E_i$ and $F = \bigotimes_{i=1}^r F_i$, where $E_i$ and $F_i$ are simple $\fg^{M_i}$-modules for all $i$.

Let $f_E$ and $f_F$ be the finitely supported $\Gamma$-invariant sections from $\Max S$ to $\mathcal{R}$ corresponding to the isomorphism classes of $E$ and $F$, respectively. Then we have the following description of extensions:

\begin{itemize}
\item[\em(1)] If $f_E$ and $f_F$ differ on more than one $\Gamma$-orbit of $\Max S$, then
$$ \Ext_\LL^1(E,F) =0. $$
\item[\em(2)] If $f_E$ and $f_F$ differ on the single $\Gamma$-orbit of $M_j \in \Max S$, then
$$ \Ext_\LL^1(E,F) \cong \Ext_\LL^1(E_j,F_j). $$
\item[\em(3)] If $f_E = f_F$, then
$$ \big((\LL/\LL')^*\big)^{\oplus\,r-1} \oplus \Ext_\LL^1(E,F) \cong \bigoplus_{i=1}^r \Ext_\LL^1(E_i,F_i). $$
\end{itemize}
\end{theorem}\qed
%\proof The proof of \cite[Theorem 3.7]{NS} generalizes to our setting in the obvious way.\qed

\begin{remark} {\em Theorem \ref{keythmext} was first proved in the simpler setting of untwisted current algebras $\fg \otimes_k S$ in \cite{Kodera}.}
\end{remark}

Theorem \ref{keythmext} reduces the problem of understanding the $\Ext_\LL^1(-,-)$ bifunctor on a pair of evaluation modules to that of understanding it on a pair of evaluation modules that are supported on a single $\Gamma$-orbit of $\Max S$. The next proposition precisely addresses this point.

\begin{proposition} \label{singleorbit} Let $k_\lambda \otimes V_f$ and $k_\mu \otimes V_g$ be simple objects of $\LL_{\ev}^{\Gamma.M}$-$\mathbf{mod}$ where $(\lambda,f)$, $(\mu,g) \in \LL^* \times \secr^\Gamma$ satisfy the three conditions of Theorem \ref{threeconditions}. In particular, $\lambda, \mu \in (\LL/\LL')_\ev^*$ , the support of $k_\lambda$, $k_\mu$,$V_f$ and $V_g$ is contained in $\Gamma.M$ and we have $V_f = V_{f(M)}$ and $V_g = V_{g(M)}$. Note that $V_{f(M)} \cong V_{g(M)} \Leftrightarrow f(M) = g(M)$ by definition. Set $U = k_{\mu-\lambda} \otimes (V_f)^* \otimes V_g$. Then there is a natural isomorphism
$$ \Ext_\LL^1(k_\lambda \otimes V_f,k_\mu \otimes V_g) \cong
\left\{\begin{array}{ll} \Hom_{\LL/K_M}(N_M/K_M',U) & \text{if } \lambda \neq \mu \\
\Hom_{\LL/\mathfrak{Z}_M}(\tilde{N}_M/\mathfrak{Z}_M',U) & \text{if } \lambda = \mu \text{ and } f(M) \neq g(M) \\
\Hom_{\LL/\mathfrak{Z}_M}(\mathfrak{Z}_M/\mathfrak{Z}_M',U) & \text{if } \lambda = \mu \text{ and } f(M) = g(M) \end{array}\right. $$
\end{proposition}
\proof Consider the case where $\lambda \neq \mu$. The isomorphism $\Ext_\LL^1(k_\lambda \otimes V_f,k_\mu \otimes V_g) \cong \Hom_{\LL/K_M}(K_M/K_M',U)$ follows from \cite[Proposition 2.6]{NS}. Recall that $\LL/K_M \cong \fg^M$. The $\fg^M$-module $K_M/K_M'$ appears in the following short exact sequence:
$$ 0 \rightarrow N_M/K_M' \rightarrow K_M/K_M' \rightarrow K_M/N_M \rightarrow 0 $$
Induce the associated long exact sequence with the functor $\Hom_{\fg^M}(-,U)$
\begin{align}
0 \rightarrow \Hom_{\fg^M}(K_M/N_M,U&) \rightarrow \Hom_{\fg^M}(K_M/K_M',U) \rightarrow \notag \\
&\Hom_{\fg^M}(N_M/K_M',U) \rightarrow \Ext_{\fg^M}^1(K_M/N_M,U). \label{keyexactseq}
\end{align}
By Proposition \ref{structureKab} (1), $K_M/N_M$ is a trivial $\fg^M$-module. Also, since $\lambda \neq \mu$, the simple module $k_{\lambda-\mu}=k_{\mu-\lambda}^*$ is nontrivial, and it follows from Schur's Lemma that $U^{\fg^M} = 0$. Therefore,
\begin{equation} \label{firstzero}
\Hom_{\fg^M}(K_M/N_M,U) \cong 0.
\end{equation}
As $U$ is a tensor product of simple finite-dimensional modules, it is completely reducible by \cite[Theorem 4.12.1]{Etingof} and we can write $U = \bigoplus_{k = 1}^N V_k$ where the $V_k$ are finite-dimensional simple $\fg^M$-modules. Because $U^{\fg^M} = 0$, the $V_k$'s are all nontrivial. Then Proposition \ref{monpremierlemme} gives
\begin{equation} \label{secondzero}
\Ext_{\fg^M}^1(K_M/N_M,U) \cong \bigoplus_{k=1}^N \Ext_{\fg^M}^1(K_M/N_M,V_k) \cong  \bigoplus_{k=1}^N 0=0.
\end{equation}
Using the results (\ref{firstzero}) and (\ref{secondzero}) in the exact sequence (\ref{keyexactseq}), we get
$$ \Ext_\LL^1(k_\lambda \otimes V_f,k_\mu \otimes V_g) \cong \Hom_{\fg^M}(K_M/K_M',U) \cong \Hom_{\fg^M}(N_M/K_M',U). $$

For the cases where $\lambda = \mu$, $U$ is a module for $\LL/\mathfrak{Z}_M \cong [\fg^M,\fg^M]$. We use \cite[Proposition 2.6]{NS} once more to get
\begin{equation*} %\label{lastcase}
\Ext_\LL^1(k_\lambda \otimes V_f,k_\mu \otimes V_g) \cong \Hom_{[\fg^M,\fg^M]}(\mathfrak{Z}_M/\mathfrak{Z}_M',U)
\end{equation*}
When $\lambda = \mu$ and $f(M) \neq g(M)$, we finish exactly as above using the long exact sequence obtained from the functor $\Hom_{[\fg^M,\fg^M]}(-,U)$ and the short exact sequence $0 \rightarrow \tilde{N}_M/\mathfrak{Z}_M' \rightarrow \mathfrak{Z}_M/\mathfrak{Z}_M' \rightarrow \mathfrak{Z}_M/\tilde{N}_M \rightarrow 0$. \qed

%In the last case, where $\lambda = \mu$ and $f(M) = g(M)$, we cannot simplify further since $(V_{f=g}^* \otimes V_{f=g})^{[\fg^M,\fg^M]} = k_0 \neq 0$ so we keep (\ref{lastcase}) as it is.
%\qed

\bigskip

With Theorem \ref{keythmext} and Proposition \ref{singleorbit} in hand, we are ready to give a description of the extension blocks of $\LL_\ev$-$\mathbf{mod}$. For this, we start by making a few simple observations about extensions of evaluation modules and introduce a few notions.

\begin{definition} {\em Let $M \in \Max S$. Then $\Bloc(\LL \, | \, \Gamma.M)$ will denote the extension blocks of the full subcategory $\LL_\ev^{\Gamma.M}$-$\mathbf{mod}$ of $\LL_\ev$-$\mathbf{mod}$. }
\end{definition}

\begin{remark} {\em In general, the categories $\Bloc(\LL \, | \, \Gamma.M)$ and $\Bloc(\fg^M)$ are not the same. However, there is always a well-defined map
\begin{align}
&\Bloc(\fg^M) \quad \longrightarrow \quad \Bloc(\LL \, | \, \Gamma.M) \\
&\qquad \quad \llbracket V \rrbracket \qquad \longmapsto \qquad \llbracket V \rrbracket \notag
\end{align}}
\end{remark}

Analogous to the framework used in Section 2 to describe isomorphism classes of simple objects, we consider the trivial fiber bundle
$$ \B \; = \bigsqcup_{M \in \Max S} \Bloc(\LL \, | \, \Gamma.M) \twoheadrightarrow \Max S $$
%= \{ \text{ sections of } \B \text{ }\}$.
\begin{definition} {\em Define the {\em support of a section} $\chi$ of $\B$ by
$$ \supp \chi = \big\{M \in \Max S \; | \; \chi(M) \neq \llbracket k_0 \rrbracket \in \Bloc(\LL \, | \, \Gamma.M)\big\}$$}
\end{definition}

Let $M \in \Max S$ and $\gamma \in \Gamma$. As $\Gamma.M = \Gamma.{}^\gamma M$, we have the equivalences of categories \eqref{catequiv}, and $\Bloc(\LL \, | \, \Gamma.M) = \Bloc\left(\LL \, | \, \Gamma.{}^\gamma M\right)$. It follows that there is a group action $\chi \mapsto {}^\gamma \chi$ of $\Gamma$ on the set of sections of $\B$ given by the following definition:
\begin{equation} \label{actioncarspec}
{}^\gamma \chi : \; M \longmapsto \left\llbracket\Big(\rho_{{}^{\gamma^{-1}}M} \circ \big(\gamma(M)\big)^{-1},V_{{}^{\gamma^{-1}}M}\Big)\right\rrbracket.
\end{equation}
where $(\rho_{{}^{\gamma^{-1}}M},V_{{}^{\gamma^{-1}}M})$ be a simple module in the extension block $\chi({}^{\gamma^{-1}} M)$. This is well-defined, since the equivalence of categories \eqref{catequiv} preserves nontrivial extensions. Henceforth, the shorter notation $\chi({}^{\gamma^{-1}}M) \circ \big(\gamma(M)\big)^{-1}$ will be used for  ${}^\gamma \chi(M)$ instead of \eqref{actioncarspec}.

%\begin{eqnarray}\label{defactionspecchar}
%\LL_\ev^{\Gamma.M}\text{-}\mathbf{mod} & \cong & \LL_\ev^{\Gamma.({}^\gamma M)}\text{-}\mathbf{mod } =  \LL_\ev^{\Gamma.M}\text{-}\mathbf{mod} \\
%(\rho,V) & \longmapsto & \left(\rho \circ \big(\gamma(M)\big)^{-1},V\right) \notag
%\end{eqnarray}
%there is always a simple evaluation module with support equal to $\Gamma.M$ in $\chi(M)$. Fix $M$ and call such a module $V_M$. We can then see

%Let $(\rho_{{}^{\gamma^{-1}}M},V_{{}^{\gamma^{-1}}M})$ be a simple module in the extension block $\chi({}^\gamma M)$.  We then define ${}^\gamma \chi$ as the section of $\B$ given by
%\begin{equation} \label{actioncarspec}
%{}^\gamma \chi : \; M \longmapsto \left\llbracket\Big(\rho_{{}^{\gamma^{-1}}M} \circ \big(\gamma(M)\big)^{-1},V_{{}^{\gamma^{-1}}M}\Big)\right\rrbracket.
%\end{equation}

\begin{definition} {\em The {\em spectral characters} of $\LL$ are defined to be the finitely supported $\Gamma$-invariant sections of $\B$ under the action of $\Gamma$ on sections described by (\ref{actioncarspec}). We use the following notation:
$$ \carspec^\Gamma = \big\{ \text{ $\Gamma$-invariant sections of } \B \text{ with finite support } \big\}.$$}
\end{definition}

Next, we associate a spectral character to each simple object of $\LL_\ev$-$\mathbf{mod}$.  By the classification recalled in Section 2, they are in natural bijection with $\F^\Gamma$, the $\Gamma$-invariant finitely supported sections of the fibration $\R$ from line (\ref{bundleiso}).

If $V$ is a simple object of $\LL_\ev$-$\mathbf{mod}$ and $f_V \in \F^\Gamma$ is its associated $\Gamma$-invariant section, we define the spectral character of $V$ to be $\chi_V \in \carspec^\Gamma$, defined by
\begin{equation} \label{speccharsimplemod}
\chi_V(M) = \llbracket f_V(M) \rrbracket \in \Bloc(\LL\; | \; \Gamma.M),
\end{equation}
where $\llbracket f_V(M) \rrbracket$ denotes the extension block that contains the class $f_V(M) \in \Rep(\fg^M)$, where $\fg^M$-modules are reinterpreted as evaluation $\LL$-modules in the obvious way.

As in the case of loop algebras and equivariant map algebras \cite{CM,Se,NS}, spectral characters can be used to classify extension blocks of evaluation modules for twisted current algebras.

%This definition of a spectral character $\chi_V \in \carspec^\Gamma$ associated to a simple object $V$ of $\LL_\ev$-$\mathbf{mod}$ will correspond to its extension block as an $\LL$-module :

\begin{proposition} \label{extblocev} Let $V$ and $W$ be simple objects of $\LL_\ev$-$\mathbf{mod}$. Then $\llbracket V \rrbracket = \llbracket W \rrbracket$ in $\LL_{\ev}$-$\mathbf{mod}$ if and only if $\; \chi_V = \chi_W \in \carspec^\Gamma$.
\end{proposition}
%\proof Just as in \cite{NS}. \qed
%\bigskip

The next corollary is an immediate consequence of Propositions \ref{blockdef} and \ref{extblocev}.

\begin{corollary} \label{blocev} Let $\LL$ be a twisted current algebra. The extension blocks of $\LL_\ev$-$\mathbf{mod}$ are in natural bijection with the elements of $\carspec^\Gamma$. The block decomposition of $\LL_{\ev}$-$\mathbf{mod}$ is then
$$ \text{$\LL_{\ev}$-$\mathbf{mod}$} = \bigoplus_{\chi \in \carspec^\Gamma} \left(\LL_{\ev}\text{-}\mathbf{mod}\right)^\chi $$
where $\left(\LL_{\ev}\text{-}\mathbf{mod}\right)^\chi$ is the full subcategory consisting of modules whose irreducible subquotients all have spectral character $\chi$.
\end{corollary}

The results of Sections 3 and 4 now determine the block decomposition of $\LL$-$\mathbf{mod}$:

\begin{theorem} \label{blockdecomp} Let $\LL$ be a twisted current algebra. The extension blocks of $\LL$-$\mathbf{mod}$ are in natural bijection with the elements of $\carspec^\Gamma\times (\LL/\LL')_{\nev}^* \;$. The block decomposition of $\LL$-$\mathbf{mod}$ is
$$ \LL\text{-}\mathbf{mod} = \bigoplus_{(\chi,\alpha) \,\in\, \mathsmaller{\carspec}^\Gamma\times \,(\LL/\LL')_{\nev}^*} (\LL\text{-}\mathbf{mod})^{(\chi,\alpha)} $$
where $\left(\LL_{\ev}\text{-}\mathbf{mod}\right)^{(\chi,\alpha)}$ is the full subcategory consisting of modules whose irreducible subquotients all have non-evaluation part $k_\alpha$ and whose evaluation part has spectral character $\chi$.
\end{theorem}
\proof The description of extension blocks follows directly from Theorem \ref{extension-locality}, Proposition \ref{blocev}, and Corollary \ref{reducetoev}.  The blocks are then determined by the discussion in Definition \ref{blockdef}.\qed

\section{Specialization to Twisted Forms}

Section 4 extends known results for equivariant map algebras to the general setting of twisted current algebras, while providing previously unknown descriptions of extensions and blocks for twisted forms.  We now concentrate on the previously unexplored special case of twisted forms, and reinterpret block decompositions in terms of maps from $\Max S$ to the fundamental group $P/Q$ of the root system of $\fg$.  This generalizes results of Chari, Moura, and Senesi \cite{CM,Se} for untwisted and twisted loop algebras, the twisted forms where $S=\mathbb{C}[t,t^{-1}]$.

{\it Our approach is entirely different from that of \cite{CM,Se}, however.}  We make no use of Weyl modules and Drinfeld polynomials, and we reduce the determination of ext-groups to the explicit computation of relatively simple tensor product multiplicities, for which well known combinatorial formulas are available.  As a byproduct, our approach gives a simple formula for the dimension of the space of extensions between simple modules in the category.

%obtain more explicit results, reducing the computation of extensions to the determination of relatively simple tensor product multiplicities for which well known combinatorial formulas are available.

Throughout this section, $\fg$ will be a finite-dimensional simple $k$-Lie algebra, $\Gamma$ a finite group acting by automorphisms on a Galois extension $S$ of the invariant subalgebra $R=S^\Gamma$ in $k$-alg.   We fix a Cartan subalgebra $\fh\subseteq\fg$ and base $\Delta=\{\ga_1,\ldots,\ga_\ell\}$ of simple roots for the root system $Q=Q(\fg,\fh)$, and we write $\{\omega_1,\ldots,\omega_\ell\}$ and $P^+$ for the corresponding sets of fundamental weights and integral dominant weights, respectively.  The simple $\fg$-module of highest weight $\lambda\in P^+$ will be denoted by $L(\lambda)$, and $L(\lambda,M)$ will be the corresponding evaluation $\cL$-module:
$$\cL\stackrel{\ev_M}{\longrightarrow}\fg\longrightarrow\End\, L(\lambda),$$
for any maximal ideal $M$ of $S$.

\begin{theorem}\label{csillagok}
Let $\cL$ be an arbitrary $S/R$-twisted form of $\fg\ot R$.
\begin{enumerate}
\item[{\rm (1)}] Let $\lambda,\mu\in P^+$ and $M\in\Max S$.  If $\mu-\lambda$ is a simple root, then
$$\Ext_\cL^1\left(L(\lambda,M),L(\mu,M)\right)=k^{\oplus d_M},$$
where $d_M$ is the dimension of the cotangent space $I/I^2$ of $\Spec R$ at the closed point $I=M\cap R$.
\item[{\rm (2)}] Two finite-dimensional simple $\cL$-modules $E$ and $F$ are in the same ext-block if and only if the corresponding maps
$$\overline{f_E},\overline{f_F}:\ \Max S\longrightarrow P^+/\sim$$
are equal, where $\ga,\beta\in P^+$ are equivalent under $\sim$ if and only if $\ga-\beta\in Q$.  That is, the blocks of the abelian category $\cL$-mod are in natural bijection with the set of {\em spectral characters}, in the sense of \cite{CM}.
\end{enumerate}

\end{theorem}

\noindent
\proof  Let $M$ be a maximal ideal of $S$ and $I=M\cap R$.  We recall the notation of Proposition \ref{structureKab}:
\begin{align*}
K_M &= \ker (\ev_M), && I_M^\Gamma=\left(\bigcap_{\gamma\in\Gamma}{}^\gamma M\right)\cap R, \\
\mathfrak{Z}_M &= \{z \in \LL \, | \, [z,\LL] \subseteq K_M\}, && \tilde{N}_M = \{n \in \mathfrak{Z}_M \, | \, I_M^\Gamma.\,n \subseteq \mathfrak{Z}_M'\}.
\end{align*}
By \cite[Proposition 3.3]{repforms}, $K_M=I\cL$ and $\cL/K_M\cong\fg$, so $\mathfrak{Z}_M=\{z\in\cL\, |\, [z,\cL]\subseteq I\cL\}$ and $K_M$ is a maximal ideal of $\cL$.  Since every twisted form is perfect and $K_M$ is a proper ideal, it follows that $\mathfrak{Z}_M$ is a proper ideal of $\cL$ containing $K_M$.  Therefore, $\mathfrak{Z}_M=K_M=I\cL$ and $\mathfrak{Z}_M'=[\mathfrak{Z}_M,\mathfrak{Z}_M]=I^2\cL$.  By \cite[Lemmas 2.1 and 2.8]{repforms}, $\bigcap_{\gamma\in\Gamma}{}^\gamma M=IS$ and $IS\cap R=I$.  Thus $I_M^\Gamma=I$ and $\tilde{N}_M=\{n\in I\cL\, |\, I.n\subseteq I^2\cL\}$.  Clearly, $I\cL\subseteq\tilde{N}_M\subseteq I\cL$, so $\tilde{N}_M=I\cL$.

The extension $S/R$ is Galois, thus faithfully flat, and $S\ot_R\cL\cong \fg\ot_k S$ is a free $S$-module.  By \cite[Exercise 7.1]{matsumura}, this implies that $\cL$ is flat as an $R$-module.  It then follows that
$$\tilde{N}_M/\mathfrak{Z}_M'=I\cL/I^2\cL\cong (I/I^2)\ot_R\cL$$
as $R$-modules.  By the same argument, the simple Lie algebra
$$\fg\cong\cL/K_M=R\cL/I\cL\cong (R/I)\ot_R\cL$$
as $k$-Lie algebras.  The induced $\fg$-module structure on $(I/I^2)\ot_R\cL$ is then given by
$$((r+I)\ot_R z).((s+I^2)\ot_R w)=(rs+I^2)\ot_R [z,w],$$
for all $r\in R$, $s\in I$, and $z,w\in\cL$.  Reinterpreting $(I/I^2)\ot_R\cL$ as the $k$-vector space
\begin{eqnarray*}
\left(I/I^2\right)\ot_k\fg&=&\left(I/I^2\right)\ot_{R/I}\left(\left(R/I\right)\ot_R\cL\right)\\
&=&\left(\left(I/I^2\right)\ot_{R/I}\left(R/I\right)\right)\ot_R\cL\\
&=&(I/I^2)\ot_R\cL,
\end{eqnarray*}
the $\fg$-module action on $\tilde{N}_M/\mathfrak{Z}_M'$ is simply the adjoint action:
$$x.(\overline{s}\ot_k y)=\overline{s}\otimes_k[x,y],$$
for all $x,y\in\fg$ and $\overline{s}=s+I^2\in I/I^2$.  That is, as $\fg$-modules,
$$\tilde{N}_M/\mathfrak{Z}_M'\cong\fg^{\oplus d_M}.$$

Since $\cL$ is a twisted form, $\fg^M=\fg$ by \cite[Proposition 3.3]{repforms} and the proof of \cite[Theorem 3.2]{twcurr}, so $Z(\fg^M)=0$ and $[\fg^M,\fg^M]=\fg$.  But $\mathfrak{Z}_M=\tilde{N}_M=I\cL$, so in light of Theorem \ref{keythmext} and Proposition \ref{singleorbit}, we need only calculate
$$\Ext_\cL^1(V,W)=\Hom_\fg\left(\fg^{\oplus d_M},V^*\ot W\right)$$
for simple evaluation modules $V$ and $W$ at the same maximal ideal $M$ to determine the block decomposition of the category $\cL-\text{mod}=\cL_\ev-\text{mod}$.  That is,
$$\Ext_\cL^1(V,W)=k^{\oplus c(V,W)d_M},$$
where $c(V,W)$ is the multiplicity of the adjoint module $\fg$ in the tensor product $V^*\ot W$ of $\fg$-modules.

These multiplicities may be computed with the PRV formula.  See \cite[Theorem 1.1]{PY}, for instance.  Explicitly, in our case, the formula simplifies to
$$c(L(\lambda),L(\mu))=c(L(\mu)^*,L(\lambda)^*)=\dim\{v\in \fg_{\mu-\lambda}\, |\, e_i^{\lambda_i+1}.v=0\text{\ for all\ }i\},$$
where $\fg_{\mu-\lambda}=\{x\in\fg\, |\, [h,x]=(\mu-\lambda)(h)x\text{\ for all\ }h\in\fh\}$, the $e_i\in\fg_{\ga_i}$ are the Chevalley-Serre generators of $\fg$, and
$$\lambda=\sum_{i=1}^\ell \lambda_i\omega_i$$
is the expression of $\lambda$ with respect to the fundamental weights $\omega_i$.  But $\mu-\lambda$ is assumed to be a simple root $\alpha_j\in\Delta$, so $\fg_{\mu-\lambda}=\fg_{\ga_j}$ is 1-dimensional and spanned by $e_j$.  By the Serre presentation,
\begin{equation}\label{etoile}
e_i^{1-a_{ji}}.e_j=0 \text{\ for all\ }i\neq j,
\end{equation}
where $\ga_j=\sum_{i=1}^\ell a_{ji}\omega_i$.  But
$$\sum_{i=1}^\ell(a_{ji}+\lambda_i)\omega_i=\ga_j+\lambda=\mu$$
is integral dominant, so $a_{ji}+\lambda_i\geq 0$ for all $i$.  That is,
$$e_i^{\lambda_i+1}.e_j=0$$
for all $i\neq j$ by (\ref{etoile}).  Obviously, $e_i.e_i=0$, so
$$c(L(\lambda),L(\mu))=1,$$
and $\Ext_\cL^1(L(\lambda,M),L(\mu,M))=\Hom_\fg(\fg^{\oplus d_M},L(\lambda)^*\ot L(\mu))$ is of dimension $d_M$, proving (1).

In particular, two finite-dimensional simple evaluation modules $L(\ga)$ and $L(\beta)$, supported on a single $\Gamma$-orbit $\Gamma.M$, are in the same ext-block if the integral dominant highest weights $\ga,\beta\in P^+$ are in the same coset of $\fh^*$ modulo the root lattice $Q$.  Conversely, if $\ga,\beta\in P^+$ and $\ga-\beta\notin Q$, then $\fg_{\beta-\ga}=0$, so $c(L(\ga),L(\beta))=0$ by the PRV formula, and $\Ext_\cL^1(L(\ga,M),L(\beta,M))=0$.  Then by Theorem \ref{keythmext} and Proposition \ref{singleorbit}, simple evaluation $\cL$-modules $E$ and $F$ belong to the same block of the category $\cL-\text{mod}$ if and only if the corresponding maps $\overline{f_E},\overline{f_F}:\ \Max\,S\longrightarrow P^+/\sim$ are equal.\qed

\begin{remark}{\em In the special case of untwisted and twisted loop algebras, the dimension of the cotangent space $I/I^2$ is always $1$, and we immediately recover the parametrisation of blocks obtained in \cite{CM,Se} from Part (1) of Theorem \ref{csillagok}.}
\end{remark}

\begin{remark}{\em
One of the most intriguing twisted forms is the Margaux algebra $\mathcal{M}$, a $\Cx[s^{\pm 1},t^{\pm 1}]/\Cx[s^{\pm 2},t^{\pm 2}]$-twisted form of $\mathfrak{sl}_2(\Cx)\ot_\Cx\Cx[s^{\pm 2},t^{\pm 2}]$, which is not a twisted multiloop algebra.  Its existence is only known cohomologically \cite{GP}, though its finite-dimensional simple modules were classified up to isomorphism in \cite{repforms}.  In particular, they are of the form
$$L(\lambda_1,M_1)\ot\cdots\ot L(\lambda_m,M_m),$$
where $\lambda_1,\ldots,\lambda_m\in P^+=\mathbb{N}$, and the ideals $M_1,\ldots,M_m\in\Max\,\Cx[s^{\pm 1},t^{\pm 1}]=\Cx^\times\times\Cx^\times$ belong to distinct orbits under the Galois action of $\Gamma=\mathbb{Z}_2\times\mathbb{Z}_2$.  Identifying each $M_i$ with a pair of nonzero complex numbers $(a_i,b_i)$, the $\Gamma$-orbit of $M_i$ is $\{(\pm a_i,\pm b_i)\}$.

As $P/Q=\mathbb{Z}_2$ for $\mathfrak{sl}_2(\Cx)$, Theorem \ref{csillagok} gives an explicit description of the blocks in the category $\mathcal{M}$-mod: they are in bijection with the finite subsets of
$$\Cx_+\times\Cx_+\times\{0,1\},$$
where $\Cx_+\times\Cx_+=\{(z_1,z_2)\in \Cx\times\Cx\, |\, \text{im}\,z_j>0\text{\ or\ }z_j\in \mathbb{R}_{>0}\text{\ for\ }j=1,2\}$ is a set of representatives of the orbit space $\Gamma\setminus(\Cx^\times\times\Cx^\times)$ and $\{0,1\}=P^+/\sim$.}
\end{remark}

\appendix
\section{Appendix: Basic Results on Extensions}
In this appendix, we recall some well-known general facts about extensions.

%The 1-dimensional representations of a Lie algebra $L$ are in natural bijection with the dual space $(L/L')^*$.  For $\lambda \in (L/L')^*$, we denote by $k_\lambda$ its associated 1-dimensional representation. In particular, $k_0$ is the trivial 1-dimensional representation, $k_\lambda \otimes k_\mu \cong k_{\lambda + \mu}$, and $(k_\lambda)^* \cong k_{-\lambda}$, for all $\lambda,\mu \in (L/L')^*$.  %Writing $V^*=\Hom_k(V,k_0)$ for the dual of any $L$-module $V$, we have $(k_\lambda)^* \cong k_{-\lambda}$.

An {\em extension} of an $L$-module $V$ by an $L$-module $W$ is a short exact sequence (in the category of $L$-modules) of the form
$$ 0 \longrightarrow W \longrightarrow E \longrightarrow V \longrightarrow 0. $$
Two such extensions
\begin{align*}
0 \longrightarrow W \stackrel{f_1}{\longrightarrow} E_1 \stackrel{g_1}{\longrightarrow} V \longrightarrow 0\\
0 \longrightarrow W \stackrel{f_2}{\longrightarrow} E_2 \stackrel{g_2}{\longrightarrow} V \longrightarrow 0
\end{align*}
are {\em equivalent} if there is an $L$-module isomorphism $\phi:\ E_1\rightarrow E_2$ such that $f_2=\phi\circ f_1$ and $g_1=g_2\circ\phi$.  We denote by $\Ext^1_L(V,W)$ the set of equivalence classes of extensions of $V$ by $W$.

By definition, the first cohomology of a Lie algebra $L$ with coefficients in an $L$-module $V$ is the vector space $\cohomo{L}{V} = \Der(L,V) / \IDer(L,V)$, where $\Der(L,V)$ and $\IDer(L,V)$ are the spaces of derivations and inner derivations, respectively.  For $L$-modules $V$ and $W$, there is a natural vector space isomorphism between $\Ext_L^1(V,W)$ and $\gcohomo{L}{\Hom_k(V,W)}$.  See \cite{Weibel} for details. %We write $V^L=\{v\in V\ |\ x.v=0\hbox{\ for all\ }x\in L\}$ for the submodule of $L$-invariants.

\bigskip

We begin by considering extensions between $1$-dimensional representations of $L$.

\begin{lemma} \label{extdim1} Let $\lambda,\mu \in (L/L')^*$. Identify $\lambda$ and $\mu$ with the corresponding elements of $L^*$ which vanish on $L'$, and set $K = \ker(\mu - \lambda) \subseteq L$. Then there is a natural isomorphism
\begin{equation*}
\Ext^1_L(k_\lambda,k_\mu) \cong \Hom_{L/K}\left(K/K',k_{\mu - \lambda}\right).
\end{equation*}
\end{lemma}
\proof First consider the case where $\lambda = \mu$. Then $k_{\mu - \lambda} = k_0$ is trivial and $(L/L')^* \cong \cohomo{L}{k_0} \cong \Ext_L^1(k_\lambda,k_\mu)$, but in this case, $L = K$, $L/K =0$ and $(L/L')^* \cong \Hom_k(K/K',k_0)^{L/K}$ matches the formula we were aiming for.

Second, consider the case where $\lambda \neq \mu$. The set $\IDer(L,k_{\mu - \lambda})$ consists of derivations $\ell \mapsto \ell.a$ where $a \in k_{\mu - \lambda}$. We find $\IDer(L,k_{\mu - \lambda}) = \Span_k \{\mu - \lambda\}$.

On the other hand, there exists $u \in L$ such that $\big(\mu - \lambda\big)(u) = 1$ and thus, an arbitrary derivation $d \in \Der(L,k_{\mu-\lambda})$ can be uniquely written as
$$ \Big(d - d(u)\big(\mu - \lambda\big)\Big) + d(u)\big(\mu - \lambda\big). $$

This means that in $\Ext^1_L(k_\lambda,k_\mu) \cong \cohomo{L}{k_{\mu - \lambda}}$, the class of $d$ is uniquely represented by $d - d(u)\big(\mu - \lambda\big)\in\Der(L,k_{\mu-\lambda})$, a derivation that vanishes on $\Span_k\{u\}$. Therefore, as a set, $\Ext^1_L(k_\lambda,k_\mu)$ is the same as $\mathcal{D}_0 = \{d \in \Der(L,k_{\mu - \lambda}) \; | \; d(u) = 0\}$.

Fix $d \in \mathcal{D}_0$. As $K \subseteq L$ is a $k$-ideal of codimension one and since $u \notin K$, we have $L = \Span_k\{u\} \oplus K$. This decomposition for $L$ implies that $d$ is determined by its restriction $d|_K \in \Der(K,k_{\mu - \lambda})$. Then as $k_{\mu - \lambda}$ is a trivial $K$-module, $ \Der(K,k_{\mu - \lambda}) \cong \Hom_k(K/K', k_{\mu - \lambda})$. We deduce that $\mathcal{D}_0 \hookrightarrow \Hom_k\big(K / K',k_{\mu - \lambda}\big)$.

Moreover, since $d \in \mathcal{D}_0$, we have $d(\omega) = 1 \cdot d(\omega) - \omega.d(u) = d\big([u,\omega]\big)$ for all $\omega \in K$. Note that $[u,\omega] \in K$ for any given $\omega \in K$. It follows that
\begin{equation} \label{conditionext1dim}
%\mathcal{D}_0 \quad \stackrel{1:1}{\longleftrightarrow} \quad
\mathcal{D}_0 \cong \Big\{f \in \Hom_k\big(K / K',k_{\mu - \lambda}\big) \; \big| \; f\big(\omega - [u,\omega] + K'\big) = 0 \text{ for all } \omega \in K \Big\}.
\end{equation}

But the set in (\ref{conditionext1dim}) is precisely $\Hom_{L/K}\big(K / K',k_{\mu - \lambda}\big)$, since $\big(\mu - \lambda\big)(u)=1$ and $L = \Span_k\{u\} \oplus K$. \qed
%
%The last remaining step justifies that for a function $f \in \Hom_k\big(K / K',k_{\mu - \lambda}\big)$, the condition $f\big(\omega - [u,\omega] + K'\big) = 0$ is equivalent to that of $f$ being invariant under the action of $L/K$ on the $L/K$-module $\Hom_k\big(K / K',k_{\mu - \lambda}\big)$.

%Let $f \in \big(K / K'\big)^*$ that satisfies the condition from line (\ref{conditionext1dim}). Then
%\begin{align*}
%&\qquad\qquad\:\: f\big(\omega - [u,\omega] + K'\big) = 0 \text{ pour tout } \omega \in K \\
%&\Longleftrightarrow \qquad f\big(A\omega - A[u,\omega] + K'\big) = 0 \quad \text{for all } \omega \in K \text{ and } A \in k \\
%&\Longleftrightarrow \qquad f\big(A\lambda(u)\omega - [Au,\omega] + K'\big) = 0 \quad \text{for all } \omega \in K \text{ and } A \in k \\
%&\Longleftrightarrow \qquad f\Big(\big(\mu - \lambda\big)(Au + \eta)\omega - [Au + \eta,\omega] + K'\Big) = 0 \quad \text{for all } \omega, \eta \in K \text{ and } A \in k \\
%&\Longleftrightarrow \qquad \big(\mu - \lambda\big)(Au + \eta) \cdot f(\omega + K'\big) - f\big([Au + \eta,\omega] + K'\big) = 0 \quad \text{for all } \omega, \eta \in K \text{ and } A \in k \\
%&\Longleftrightarrow \qquad \big(\mu - \lambda\big)(\ell) \cdot f(\omega + K'\big) - f\big([Au + \eta,\omega] + K'\big) = 0 \quad \text{for all } \omega \in K \text{ and } \ell \in \LL
%\end{align*}

\begin{corollary}{\bf \cite[Corollary 2.5]{NS}} \label{extab} Let $k_\lambda$, $k_\mu$ be 1-dimensional modules over an abelian Lie algebra $Z$.  Then
$$ \Ext_Z^1(k_\lambda,k_\mu) \cong \left\{\begin{array}{cl} Z^* & \text{if } \lambda = \mu \\ 0 & \text{if } \lambda \neq \mu.\end{array}\right.$$

\noindent
\proof {\em The case $\lambda = \mu$ follows directly from Lemma \ref{extdim1}. However, if $\lambda \neq \mu$, the same lemma gives us $\Ext_Z^1(k_\lambda,k_\mu) \cong \Hom_{Z/K}(K/K',k_{\mu - \lambda})$ where $K = \ker (\mu - \lambda)$.

Given that $Z$ is abelian, $K/K'$ is a trivial $Z/K$-module. Then for any choice of \linebreak $f \in \Hom_{Z/K}(K/K',k_{\mu - \lambda})$, we get $f(K/K') \subseteq (k_{\mu-\lambda})^{Z/K}$. Since $k_{\mu - \lambda}$ is a nontrivial simple $Z/K$-module, $(k_{\mu-\lambda})^{Z/K} = 0$. It follows that $f = 0$. \qed}
\end{corollary}

\bigskip

As we often work with reductive Lie algebras, it will be useful to consider $\Ext^1_{L_1 \oplus L_2}(V,W)$, when $V$ and $W$ are finite-dimensional simple modules over a Lie algebra direct sum $L_1 \oplus L_2$.  Such modules $V$ and $W$ are necessarily tensor products $M_1\otimes M_2$ of finite-dimensional simple $L_i$-modules $M_i$ for $i=1,2$ \cite{bourbaki-algebre8}.

%\begin{remark} \label{pbwrmk} As a corollary to the Poincar\'e-Birkhoff-Witt theorem, we have
%$$ U(L_1 \oplus L_2) \cong U(L_1) \otimes U(L_2) $$
%From such an isomorphism, it can be deduced that every finite-dimensional simple \linebreak $(L_1 \oplus L_2)$-module is of the form $M_1 \otimes M_2$ where $M_i$ is a simple $L_i$-module. Conversely, for any pair $M_1$, $M_2$ where $M_i$ are finite-dimensional simple $L_i$-modules, the module $M_1 \otimes M_2$ is a simple $(L_1 \oplus L_2)$-module. For a justification of these results, see ETINGOF.
%\end{remark}

The following general fact is a straightforward consequence of the K\"unneth formula and the usual duality between homology and cohomology.  See \cite[Proposition 2.7b]{NS} for details.

%is tous ces résultats sur les extensions entre modules de dimension 1, voici un résultat cohomologique relativement général, mais bien utile :

\begin{proposition} \label{extsommelie} Let $L_1$ and $L_2$ be Lie algebras and let $V = V_1 \otimes V_2$ and \linebreak $W = W_1 \otimes W_2$ be finite-dimensional simple $(L_1 \oplus L_2)$-modules. Then
\begin{equation*}
\Ext_{L_1 \oplus L_2}^1(V, W) \cong \left\{
\begin{array}{cl} 0 & \text{if } V_1 \ncong W_1 \text{ and } V_2 \ncong W_2\\
\Ext_{L_2}^1(V_2,W_2) & \text{if } V_1 \cong W_1 \text{ and } V_2 \ncong W_2\\
\Ext_{L_1}^1(V_1,W_1) & \text{if } V_1 \ncong W_1 \text{ and } V_2 \cong W_2\\
\Ext_{L_1}^1(V_1,W_1) \oplus \Ext_{L_2}^1(V_2,W_2) & \text{if } V_1 \cong W_1 \text{ and } V_2 \cong W_2.
\end{array} \right.
\end{equation*}
%\proof See proposition 2.7 b) in \cite{NS}. The proof uses the so called K\"unneth formula for homology and the duality notions that link Lie algebra homology and cohomology. \qed
\end{proposition}

The following result is needed for the application to twisted current algebras in Section 4.

\begin{proposition} \label{monpremierlemme} Let $T$ and $V$ be modules over a finite-dimensional reductive $k$-Lie algebra $\fg$, where $T$ is trivial and $V$ is finite-dimensional, nontrivial, and simple.  Then
$$ \Ext_\fg^1(T,V) = 0. $$
\end{proposition}
\proof Since $\Ext_\fg^1(T,V) \cong \gcohomo{\fg}{\Hom_k(T,V)}$, it suffices to show that all derivations in $\Der\big(\fg,\Hom_k(T,V)\big)$ are inner.

Let $d \in \Der\big(\fg,\Hom_k(T,V)\big)$. For each $x,y\in \fg$ and $t\in T$, we have
$$\Big(d\big([x,y]\big)\Big)(t)=x.\Big(\big(d(y)\big)(t)\Big) - y.\Big(\big(d(x)\big)(t)\Big),$$
so each $t\in T$ induces a derivation $\delta_t \in \Der(\fg,V)$ by setting $\delta_t(x)=\big(d(x)\big)(t)$.

%$d \in \Hom_k\big(R,\Hom_k(T,V)\big)$ and for all $x,y \in L$ and all $t \in T$, we have
%\begin{align*}
%\Big(d\big([x,y]\big)\Big)(t) &= \big(x \cdot d(y)\big)(t) - \big(y \cdot d(x)\big)(t)\\
%&= x.\Big(\big(d(y)\big)(t)\Big) - \big(d(y)\big)(x.t) + y.\Big(\big(d(x)\big)(t)\Big) - \big(d(x)\big)(y.t)\\
%&= x.\Big(\big(d(y)\big)(t)\Big) - y.\Big(\big(d(x)\big)(t)\Big)
%\end{align*}
%Fix $t \in T$. The above defines a derivation of $\delta_t \in \Der(R,V)$ given by $x \mapsto \big(d(x)\big)(t)$.

Next, we will prove that $\delta_t$ is inner.  Recall that
$$\cohomo{\fg}{V}=\cohomo{\fg}{k_0^*\ot V}=\gcohomo{\fg}{\Hom_k(k_0,V)}=\Ext_\fg^1(k_0,V),$$
so it suffices to show that $\Ext_\fg^1(k_0,V) = 0$. %As a first step in doing so, let's mention that since the Lie algebra $\fg$ is finite-dimensional and reductive over an algebraically closed field of characteristic zero,
Since $\fg$ is reductive,
$$ \fg = \mathcal{Z} \oplus \mathcal{S}, $$
where $\mathcal{Z}= Z(\fg)$ and $\mathcal{S} = \fg'$ is semisimple (and finite-dimensional).  We can thus decompose $V$ as $V = k_\lambda \otimes V_\mathcal{S}$ where $\lambda \in \mathcal{Z}^*$ and $V_\mathcal{S}$ is a finite-dimensional simple $\mathcal{S}$-module.

We now apply Proposition \ref{extsommelie} with $L_1 = \mathcal{Z}$ and $L_2 = \mathcal{S}$:
\begin{equation} \label{omgpourlewin}
\Ext_{\fg}^1(k_0, V) \cong \left\{
\begin{array}{cl} 0 & \text{if } 0 \neq \lambda \text{ and } k_0 \ncong V_\mathcal{S} \\
\Ext_\mathcal{S}^1(k_0,V_\mathcal{S}) & \text{if } 0 = \lambda \text{ and } k_0 \ncong V_\mathcal{S} \\
\Ext_\mathcal{Z}^1(k_0,k_\lambda) & \text{if } 0 \neq \lambda \text{ and } k_0 \cong V_\mathcal{S} \\
\Ext_\mathcal{Z}^1(k_0,k_\lambda) \oplus \Ext_\mathcal{S}^1(k_0,V_\mathcal{S}) & \text{if } 0 = \lambda \text{ and } k_0 \cong V_\mathcal{S}
\end{array} \right.
\end{equation}
By Weyl's Theorem, $\Ext_\mathcal{S}^1(k_0,V_\mathcal{S}) =0$.
%On the other hand, Corollary \ref{extab} gives values for the sets $\Ext_\mathcal{Z}^1(k_0,k_\lambda)$.
Combining Corollary \ref{extab} with (\ref{omgpourlewin}), we obtain
\begin{equation} \label{omgpourlewinv2}
\Ext_{\fg}^1(k_0, V) \cong \left\{
\begin{array}{cl} 0 & \text{if } V \text{ is a nontrivial } \fg\text{-module} \\
\mathcal{Z}^* & \text{if } V \text{ is a trivial } \fg\text{-module}.
\end{array} \right.
\end{equation}
But since $V$ is nontrivial by assumption, we see that $\Ext_{\fg}^1(k_0, V) =0$, so each derivation $\delta_t$ is inner.
%. It means that for any, given $t \in T$ the derivation $\delta_t$ belongs to the set $\IDer(\fg,V)$ of interior derivations.

Next, by fixing a vector space section $s:\ \IDer(\fg,V)\rightarrow V$ of the short exact sequence
$$ 0\rightarrow V^\fg\rightarrow V\rightarrow\IDer(\fg,V)\rightarrow 0,$$
we have $\delta_t(x) =x.s(\delta_t)$ for all  $x \in \fg$.  The set-theoretic mapping $f_d : t \mapsto \delta_t \mapsto s(\delta_t)$ is clearly $k$-linear. %Therefore, we can write $f_d \in \Hom_k(T,V)$.

Finally, since $T$ is trivial,
\begin{align*}
\big(d(x)\big)(t) &= \delta_t(x) = x.f_d(t) = x.f_d(t) - f_d(x.t) = \big(x \cdot f_d\big)(t),
\end{align*}
for all $x\in\fg$ and $t\in T$.  Thus $d(x) = x \cdot f_d$ for all $x\in\fg$, and every derivation in $\Der\big(\fg,\Hom_k(T,V)\big)$ is inner. \qed

%$d \in \IDer\big(\fg,\Hom_k(T,V)\big)$ i.e. $d$ is an interior derivation and this proof is complete. \qed

\end{document}